\documentclass[11pt]{article}  
\usepackage[margin=1in]{geometry}
\usepackage{graphicx}
\usepackage{url}
\usepackage[utf8]{inputenc}
\usepackage[english]{babel}
 
\usepackage{amsthm}
\usepackage{etoolbox}
\usepackage{comment}
\usepackage{amsmath}
\usepackage[utf8]{inputenc}
\usepackage{amssymb}
\usepackage{xcolor}
\usepackage{enumerate}
\usepackage{mathtools}
\usepackage{graphicx}

\newcommand{\beq}{\begin{equation}}
\newcommand{\eeq}{\end{equation}}
\newcommand{\beqa}{\begin{eqnarray}}
\newcommand{\eeqa}{\end{eqnarray}}

\newtheorem{remark}{\bf Remark}[section]

\title{From a discrete model of chemotaxis with volume-filling to a generalised Patlak-Keller-Segel model}

\author{
Federica Bubba$^{1}$, Tommaso Lorenzi$^{2,3}$ and Fiona R. Macfarlane$^{2}$}

\date{$^{1}$Sorbonne Universit\'es, Universit\'es Paris-Diderot, Laboratoire Jacques-Louis Lions, F-75005 Paris, France\\
$^{2}$School of Mathematics and Statistics, University of St Andrews, St Andrews, KY16 9SS, United Kingdom\\
$^{3}$ Department of Mathematical Sciences ``G. L. Lagrange'', Dipartimento di Eccellenza 2018-2022, Politecnico di Torino, 10129 Torino, Italy}

\begin{document}


\maketitle
\begin{abstract}
We present a discrete model of chemotaxis whereby cells responding to a chemoattractant are seen as individual agents whose movement is described through a set of rules that result in a biased random walk. In order to take into account possible alterations in cellular motility observed at high cell densities (\emph{i.e.} volume-filling), we let the probabilities of cell movement be modulated by a decaying function of the cell density. We formally show that a general form of the celebrated Patlak-Keller-Segel (PKS) model of chemotaxis can be formally derived as the appropriate continuum limit of this discrete model. The family of steady-state solutions of such a generalised PKS model are characterised and the conditions for the emergence of spatial patterns are studied via linear stability analysis. Moreover, we carry out a systematic quantitative comparison between numerical simulations of the discrete model and numerical solutions of the corresponding PKS model, both in one and in two spatial dimensions. The results obtained indicate that there is excellent quantitative agreement between the spatial patterns produced by the two models. Finally, we numerically show that the outcomes of the two models faithfully replicate those of the classical PKS model in a suitable asymptotic regime.
\end{abstract}

\section{Introduction}
The ability of living organisms to form self-organised spatial patterns is at the root of a wide range of developmental and evolutionary phenomena~\cite{COUZIN20031,johnson2010self}. In many biological systems, the emergence of spatial organisation is orchestrated by chemotaxis, whereby individuals undergo directed migration in response to the gradient of chemical signals (\emph{i.e.} chemoattractants)~\cite{van2004chemotaxis}. Chemotaxis plays a pivotal role in many biological processes -- such as the immune response to infection, wound healing, embryogenesis, cancer progression and metastasis~\cite{constantin2010leukocyte,kay2008changing,kundra1994regulation,raina2019role,roussos2011chemotaxis,wadhams2004making} -- and the mathematical modelling of chemotactic movement has received considerable attention from mathematicians and physicists over the past seventy years -- the interested reader is referred to~\cite{chalub2006model,hillen2009user,painter2019mathematical,perthame2006transport} and references therein.

\paragraph{The Patlak-Keller-Segel model of chemotaxis} The simplest and most classical mathematical model for the emergence of self-organised spatial patterns driven by chemotaxis in biological systems (\emph{e.g.} populations of bacteria and eukaryotic cells) dates back to Patlak~\cite{Patlak1953} and Keller-Segel~\cite{Keller1970}. This model comprises a conservation equation for the density of cells and a balance equation for the concentration of chemoattractant in the form of the following system of coupled parabolic equations
\begin{equation}
		\begin{cases}
\displaystyle{\frac{\partial u}{\partial t} - \nabla \cdot \left(\beta_u \nabla u - \chi \, u \, \nabla c  \right) = 0, \quad u \equiv u(t,x),}
\\\\
\displaystyle{\frac{\partial c}{\partial t} - \beta_c \, \Delta c = \alpha \, u - \kappa \, c, \quad\quad\quad\quad\;\; c \equiv c(t,x),}
		\end{cases}
(t,x) \in \mathbb{R}^+_{*} \times \Omega,	
		\label{eq:KellerSegelSystemOri}
\end{equation}
subject to biologically relevant initial and boundary conditions. Here, the real, non-negative functions $u(t,x)$ and $c(t,x)$ represent, respectively, the density of cells and the concentration of chemoattractant at time $t \in \mathbb{R}^+_{*}$ and at position $x \in \Omega$. The set $\Omega$ is an open and bounded subset of $\mathbb{R}^d$ with smooth boundary $\partial \Omega$ and $d=1,2,3$ depending on the biological problem under study.

In the Patlak-Keller-Segel (PKS) model~\eqref{eq:KellerSegelSystemOri}, the transport term in the equation for $u$ models the rate of change of the cell density due to both undirected, random movement and chemotaxis. Undirected, random movement is described through Fick's first law of diffusion with diffusivity $\beta_u>0$. Furthermore, chemotaxis is modelled via an advection term whereby the velocity field is proportional to $\nabla c$, in order to capture the tendency of cells to move toward regions of higher concentration of the chemoattractant (\emph{i.e.} cells move up the gradient of the chemoattractant). The proportionality constant $\chi > 0$ represents the sensitivity of cells to the chemoattractant (\emph{i.e.} the chemotactic sensitivity). The second term of the left-hand side of the equation for $c$ models the effect of Fickian diffusion and $\beta_c>0$ is the diffusivity of the chemoattractant. Moreover, the first term on the right-hand side takes into account the fact that the chemoattractant is produced by the cells at a rate $\alpha > 0$, while the second term describes natural decay of the chemoattractant, which occurs at rate $\kappa >0$.  

Although the PKS model~\eqref{eq:KellerSegelSystemOri} has helped to elucidate the mechanisms that underlie the formation of self-organised spatial patterns in various biological contexts~\cite{painter2019mathematical}, it is well known that its solutions may blow up in finite time~\cite{calvez2008parabolic, nagai2001global, Winkler2013}. In order to avoid unphysical finite-time blow-up, a number of possible variations on the classical PKS model have been proposed in the literature~\cite{painter2019mathematical,hillen2009user}. In a number of these variants, the cell diffusivity and sensitivity to the chemoattractant are assumed to be functions of the cell density, leading to modified versions of~\eqref{eq:KellerSegelSystemOri} of the following form
\begin{equation}
		\begin{cases}
\displaystyle{\frac{\partial u}{\partial t} - \nabla \cdot \left(\beta_u \, D(u) \, \nabla u -  \chi \, \psi(u) \, u \, \nabla c  \right) = 0, \quad u \equiv u(t,x),}
\\\\
\displaystyle{\frac{\partial c}{\partial t} - \beta_c \, \Delta c = \alpha \, u - \kappa \, c, \quad\quad\quad\quad\quad\quad\quad\quad\quad\;\; c \equiv c(t,x),}
		\end{cases}
(t,x) \in \mathbb{R}^+_* \times \Omega.		
		\label{eq:KellerSegelSystem}
\end{equation}
Compared to the classical PKS model~\eqref{eq:KellerSegelSystemOri}, here undirected, random cell movement is modelled as a nonlinear diffusion process with diffusivity $\beta_u \, D(u)$ and the chemotactic sensitivity is a function of the cell density $\chi \, \psi(u)$. The solutions are prevented from blowing up in finite time by introducing suitable assumptions on the functions $D(u)$ and $\psi(u)$. In particular, Hillen and Painter~\cite{Painter2002} have shown that letting
\begin{equation}
\label{ass:volfil}
D(u) := \psi(u) - u \, \psi'(u)
\end{equation}
and assuming $\psi(u)$ to be a monotonically decreasing function of the cell density (\emph{i.e.} $\psi'(\cdot) \leq 0$), enables to capture possible alterations in cellular motility observed at high cell densities (\emph{i.e.} volume-filling) and precludes blow-up from occurring. Moreover, under the additional assumptions
\begin{equation*}
\psi(u) > 0 \; \text{ for } \; 0 \leq u < \overline{u} \quad \text{and} \quad \psi(\overline{u}) = 0,
\end{equation*}
where $\bar{u}>0$ is a critical value of the cell density above which no more cells can move into a given position, the same authors have proven global existence of classical solutions of~\eqref{eq:KellerSegelSystem} subject to suitable initial and boundary conditions~\cite{Hillen2001}. More recently, focussing on the case where $D(u)$ is a constant function and building upon the modelling strategy presented by Painter \emph{et al.}~\cite{painter2018}, Bubba \emph{et al.}~\cite{bubba2019chemotaxis} considered the following definition of the sensitivity of cells to the chemoattractant
\begin{equation}
\label{eq:exponential_chemofun}
	\psi(u) := \exp{(-u/u_{\text{max}})},
\end{equation}
where $u_{\text{max}}>0$ represents a critical cell density above which chemotactic movement is reduced due to overcrowding. In~\cite{painter2018,bubba2019chemotaxis} it has been shown that letting chemotactic sensitivity be an exponentially decaying function of the cell density enables to reproduce experimental results on cell pattern formation. 

\paragraph{Derivation of continuum models for the movement of living organisms from discrete models} Continuum models for the movement of living organisms like~\eqref{eq:KellerSegelSystemOri} and~\eqref{eq:KellerSegelSystem}-\eqref{eq:exponential_chemofun} are amenable to both numerical and analytical approaches, which support a more in-depth theoretical understanding of the application problems under study. However, defining these models on the basis of population-level phenomenological assumptions makes it difficult to represent fine details of the movement of single individuals. Therefore, it is desirable to derive them from first principles as the appropriate continuum limit of discrete models that track the dynamics of individual organisms. In fact, such discrete models enable a more direct and precise description of the spatial dynamics of living systems at the individual level, thus ensuring that key biological aspects are faithfully mirrored in the structure of the equations that compose the continuum model. As a consequence, the derivation of continuum models for the movement of living organisms from underlying discrete models has become an active research field. Examples in this fertile area of research include the derivation of continuum models of chemotaxis from velocity-jump and space-jump processes~\cite{othmer1988models,othmer2000diffusion,Painter2002,painter2003modelling, Charteris2014} or from different types of random walks~\cite{burger2011continuous,stevens2000derivation,stevens1997aggregation}; the derivation of diffusion and nonlinear diffusion equations from random walks~\cite{champagnat2007invasion,chaplain2019bridging,inoue1991derivation,oelschlager1989derivation,othmer2002diffusion,penington2011building,penington2014interacting}, from systems of discrete equations of motion~\cite{baker2019free,byrne2009individual,oelschlager1990large,lorenzi2019free,motsch2018short,murray2009discrete,murray2012classifying}, from discrete lattice-based exclusion processes~\cite{binder2009exclusion,dyson2012macroscopic,fernando2010nonlinear,johnston2017co,johnston2012mean,landman2011myopic,lushnikov2008macroscopic,simpson2010cell} and from cellular automata~\cite{deroulers2009modeling,drasdo2005coarse,simpson2007simulating}; and the derivation of non-local models of cell-cell adhesion from position-jump processes~\cite{buttenschoen2018space}.

\paragraph{Contents of the paper} 
In this paper, we present a discrete model of chemotaxis whereby cells responding to a chemoattractant are described as individual agents. Volume-filling effects are taken into account by modulating the probabilities of cell movement by a decaying function of the cell density, which is defined according to~\eqref{eq:exponential_chemofun} with $u_{\text{max}}$ being interpreted, in a broader biological sense, as a critical value of the cell density above which cellular movement is impaired.

In our model, cells move according to a set of rules that result in a discrete-time biased random walk on a regular lattice, which is coupled with a discrete balance equation for the concentration of chemoattractant. The modelling approach adopted here to describe cell dynamics shares some similarities with the one presented in~\cite{Charteris2014}, where also cell proliferation and cell-cell adhesion have been considered. However, the model in~\cite{Charteris2014} relies on an exclusion-based approach whereby each lattice site can be occupied by at most one cell. On the other hand, in the model considered here the maximum occupancy for a lattice site is linked to the value of the critical cell density $u_{\text{max}}$. In this regard, we formally show (see Appendix~A) that the continuum limit of our discrete model is given by~\eqref{eq:KellerSegelSystem} complemented with~\eqref{ass:volfil} and~\eqref{eq:exponential_chemofun}. Although, to the best of our knowledge, proving the global existence of solutions to~\eqref{eq:KellerSegelSystem}-\eqref{eq:exponential_chemofun} is still an open problem, the numerical solutions presented here indicate that the value of the cell density remains bounded and blow-up does not occur unless $u_{\text{max}} \to \infty$, which is the asymptotic regime in which the generalised PKS model~\eqref{eq:KellerSegelSystem}-\eqref{eq:exponential_chemofun} formally reduces to the classical PKS model~\eqref{eq:KellerSegelSystemOri}.

The paper is organised as follows. In Section~\ref{sec:hybrid}, we present our discrete model of chemotaxis with volume-filling effects. In Section~\ref{sec:comparison}, we characterise the family of the steady-state solutions to the generalised PKS model~\eqref{eq:KellerSegelSystem}-\eqref{eq:exponential_chemofun} and study the conditions for the emergence of spatial patterns via linear stability analysis. Moreover, we carry out a systematic quantitative comparison between the results of numerical simulations of the discrete model and numerical solutions of its continuum counterpart given by~\eqref{eq:KellerSegelSystem}-\eqref{eq:exponential_chemofun}. The results obtained indicate that there is excellent quantitative agreement between the spatial patterns produced by the two models, in the presence of sufficiently large cell numbers. Finally, we numerically show that the outcomes of the two models  faithfully replicate those of the classical PKS model~\eqref{eq:KellerSegelSystemOri} in the asymptotic regime whereby $u_{\text{max}} \to \infty$. Section~\ref{sec:fin} concludes the paper and provides a brief overview of possible research perspectives.

\section{From a discrete model of chemotaxis with volume-filling effects to a generalised Patlak-Keller-Segel model}
\label{sec:hybrid}	
In this section, we develop a discrete model of chemotaxis with volume-filling effects. In our model, each cell is seen as an agent that occupies a position on a lattice, while the concentration of chemoattractant is described by a discrete, non-negative function. Cells undergo undirected, random movement and chemotactic movement in the presence of volume-filling effects, while the chemoattractant is produced by the cells, undergoes natural decay and diffuses according to Fick's first law of diffusion. 

For ease of presentation, we let the cells and the chemoattractant be distributed along the real line $\mathbb{R}$, but there would be no additional difficulty in considering bounded spatial domains or higher spatial dimensions. We discretise the time variable $t \in \mathbb{R}^+$ and the space variable $x \in \mathbb{R}$ as $t_{k} = k \tau$ with $k\in\mathbb{N}_{0}$ and $x_{i} = i  h$ with $i \in \mathbb{Z}$, respectively, where $\tau, h >0$. Moreover, we introduce the dependent variable $n^{k}_{i}\in\mathbb{N}_0$ to model the number of cells on the lattice site $i$ and at the time-step $k$, and we define the corresponding density of cells as
\begin{equation}
\label{e:n}
u^k_{i} := n^{k}_{i} \,  h^{-1}.
\end{equation}
The concentration of chemoattractant on the lattice site $i$ and at the time-step $k$ is modelled by the discrete, non-negative function $c^{k}_{i}$.

\subsection{Dynamic of the chemoattractant}
We let $\beta_c >0$ be the diffusivity of the chemoattractant and we denote by $\alpha>0$ and $\kappa>0$ the rate at which the chemoattractant is produced by the cells and the rate at which it undergoes natural decay, respectively. With this notation, letting the operator $\mathcal{L}$ be the finite-difference Laplacian on the lattice $\left\{x_i\right\}_{i \in \mathbb{Z}}$, the principle of mass balance gives the following equation for the concentration of chemoattractant $c^{k}_{i}$
\begin{equation}
\label{e:cdiscrete}
c^{k+1}_{i} = c^{k}_{i} + \tau \ \left(\beta_c (\mathcal{L} \ c^k)_i + \alpha \ u^k_{i} - \kappa \ c^{k}_{i} \right),
\end{equation}
subject to zero-flux boundary conditions.

\subsection{Dynamic of the cells}
We let the cells update their positions according to a combination of undirected, random movement and chemotactic movement, which are seen as independent processes. This results in the following rules which govern the dynamic of the cells.
		
\paragraph{Mathematical modelling of chemotactic cell movement with volume-filling effects} 
Building upon the modelling strategy presented in~\cite{chaplain2019bridging}, we model chemotactic cell movement as a biased random walk whereby the movement probabilities depend on the difference between the concentration of chemoattractant at the site occupied by a cell and the concentration of chemoattractant at the neighbouring sites. Moreover, we multiply the movement probabilities by a monotonically decreasing function of the cell density at the neighbouring sites, in order to take into account volume-filling effects 
consisting in possible reduction of chemotactic sensitivity at higher cell densities. In particular, for a focal cell on the lattice site $i$ and at the time-step $k$, we define the probability of moving to the lattice site $i-1$ (\emph{i.e.} the probability of moving left) via chemotaxis as
\begin{equation}
\label{e:left}
J^{k}_{{\rm L} i} :=  \eta\ \psi(u_{i-1}^{k}) \frac{(c^{k}_{i-1}-c^{k}_{i})_{+}}{2 \, \overline{c}}, 
\end{equation}
where $(\cdot)_{+}$ denotes the positive part of $(\cdot)$, the probability of moving to the lattice site $i+1$ (\emph{i.e.} the probability of moving right) via chemotaxis as
\begin{equation}
\label{e:right}
J^{k}_{{\rm R} i} :=  \eta\ \psi(u_{i+1}^{k}) \frac{(c^{k}_{i+1}-c^{k}_{i})_{+}}{2 \, \overline{c}},
\end{equation}
and the probability of not undergoing chemotactic movement as 
\begin{equation}
\label{e:stay}
1 - J^{k}_{{\rm L} i} - J^{k}_{{\rm R} i}.
\end{equation}
Here, the weight function $\psi$ is defined according to~\eqref{eq:exponential_chemofun}, the parameter $\eta>0$ is directly proportional to the chemotactic sensitivity and we assume $\eta \, \psi(\cdot) \leq 1$. Moreover, the parameter $\overline{c}>0$ is directly proportional to the maximal value that can be attained by the concentration of chemoattractant. Dividing by $\overline{c}$ ensures that the values of the quotients in~\eqref{e:left}-\eqref{e:stay} are all between 0 and 1. In particular, the results of numerical simulations presented in Section~\ref{sec:comparison} indicate that a suitable definition of $\overline{c}$ is
\begin{equation}
\label{e:defbarc}
\overline{c} := \max \left(\max_{i \in \mathbb{Z}}  c^0_i, \zeta \ u_{{\rm max}}  \right),
\end{equation}
where $u_{{\rm max}}$ is given by~\eqref{eq:exponential_chemofun} and $\zeta>0$ is a scaling factor ensuring unit consistency. Notice that definitions~\eqref{e:left} and \eqref{e:right} are such that cells will move up the gradient of the chemoattractant.
	
\paragraph{Mathematical modelling of undirected, random cell movement with volume-filling effects} We model undirected, random cell movement as a random walk with movement probability $0~<~\theta~\leq~1$. In order to capture volume-filling effects consisting in possible reduction of cell motility at higher cell densities~\cite{gole2011quorum, d2018collective}, as similarly done in the case of chemotactic movement, we modulate the movement probability by a decreasing function of the cell density at the neighbouring sites. In particular, for a focal cell on the lattice site $i$ and at the time-step $k$, we define the probability of moving to the lattice site $i-1$ via undirected, random movement as
\begin{equation}
\label{e:leftdif}
T^{k}_{{\rm L} i} :=  \frac{\theta}{2}\ \psi(u_{i-1}^{k}),
\end{equation}
the probability of moving to the lattice site $i+1$ via undirected, random movement as	
\begin{equation}
\label{e:rightdiff}
T^{k}_{{\rm R} i} :=  \frac{\theta}{2}\ \psi(u_{i+1}^{k}),
\end{equation}
and the probability of not undergoing undirected, random movement as 
\begin{equation}
\label{e:staydiff}
1 - T^{k}_{{\rm L} i}  - T^{k}_{{\rm R} i} .
\end{equation}
In~\eqref{e:leftdif} and~\eqref{e:rightdiff}, the modulating function $\psi$  is defined according to~\eqref{eq:exponential_chemofun}.
		
\section{Comparison between discrete and continuum models}
\label{sec:comparison}
Letting $\tau, h \to 0$ in such a way that
\begin{equation}
\label{ass:quottozero1d}
\frac{\eta h^2}{2\tau \overline{c}}\rightarrow \chi \in \mathbb{R}^+_* \quad \text{and} \quad \frac{\theta h^2}{2\tau}\rightarrow \beta_{u} \in \mathbb{R}^+_* \quad \text{as} \quad \tau, h \to 0,
\end{equation}
one can formally show (see Appendix~A) that the continuum counterpart of the discrete model presented in Section~\ref{sec:hybrid} is given by the generalised PKS model~\eqref{eq:KellerSegelSystem} posed on $\mathbb{R}^+_* \times \mathbb{R}$, and complemented with~\eqref{ass:volfil} and~\eqref{eq:exponential_chemofun}. Similarly, in the case where the cells and the chemoattractant are distributed over $\mathbb{R}^2$, considering a two-dimensional regular spatial grid of step $h$, defining the cell density via the two-dimensional analogue of~\eqref{e:n}, letting the operator $\mathcal{L}$ be the finite-difference Laplacian on the grid $\left\{x_{1 i}\right\}_{i \in \mathbb{Z}} \times \left\{x_{2 j}\right\}_{j \in \mathbb{Z}}$ and assuming $\tau, h \to 0$ in such a way that
\begin{equation}
\label{ass:quottozero2d}
\frac{\eta h^2}{4\tau \overline{c}}\rightarrow \chi \in \mathbb{R}^+_* \quad \text{and} \quad \frac{\theta h^2}{4\tau}\rightarrow \beta_{u} \in \mathbb{R}^+_*, \quad \text{as } \tau, h \to 0,
\end{equation}
it is possible to formally obtain the generalised PKS model~\eqref{eq:KellerSegelSystem} posed on $\mathbb{R}^+_* \times \mathbb{R}^2$, and complemented with~\eqref{ass:volfil} and~\eqref{eq:exponential_chemofun}, as the continuum limit of our discrete model.

In this section, we carry out a systematic quantitative comparison between our discrete model and the generalised PKS model~\eqref{eq:KellerSegelSystem}-\eqref{eq:exponential_chemofun} -- \emph{i.e.} the following system of coupled parabolic equations
\begin{equation}
\label{modexpli}
\begin{cases}
\displaystyle{\frac{\partial u}{\partial t} - \nabla \cdot \left[\beta_u \, \exp{(-u/u_{\text{max}})} \, \left(1 + \frac{u}{u_{\text{max}}}\right)  \, \nabla u -  \chi \, \exp{(-u/u_{\text{max}})} \, u \, \nabla c  \right] = 0,}
\\\\
\displaystyle{\frac{\partial c}{\partial t} - \beta_c \, \Delta c = \alpha \, u - \kappa \, c,}
\end{cases}	
\end{equation}
which is obtained by substituting~\eqref{ass:volfil} and~\eqref{eq:exponential_chemofun} into~\eqref{eq:KellerSegelSystem}. In Section~\ref{sec:comparison}\ref{sec:stability}, we characterise the family of steady-state solutions of~\eqref{eq:KellerSegelSystem}-\eqref{eq:exponential_chemofun} and investigate, via linear stability analysis of the unique positive homogeneous steady state, the conditions for the emergence of spatial patterns. In Section~\ref{sec:comparison}\ref{sec:numres}, we compare the results of numerical simulations of the discrete model with numerical solutions of the continuum model given by~\eqref{eq:KellerSegelSystem}-\eqref{eq:exponential_chemofun}, both in one and in two spatial dimensions. 

\subsection{Steady-state solutions of the generalised PKS model and linear stability analysis}
\label{sec:stability}
We consider the case where~~\eqref{eq:KellerSegelSystem} is subject to an initial condition of components 
\begin{equation} \label{eq:ICs}
u(0,x) = u^0(x) \quad \text{and} \quad c(0,x) = c^0(x), \quad x \in \Omega,
\end{equation}
and to the no-flux boundary conditions
\begin{equation} \label{eq:BCs}
\nabla u(t,x) \cdot \hat{{\rm n}} = 0 \quad \text{and} \quad \nabla c(t,x) \cdot \hat{{\rm n}} = 0, \quad (t,x) \in \mathbb{R}^+_* \times \partial \Omega.
\end{equation}
Here, $u^0 \not\equiv 0$ and $c^0 \not\equiv 0$ are real, non-negative and sufficiently regular functions, and $\hat{{\rm n}}$ is the unit normal to $\partial \Omega$ that points outward from $\Omega$. 

\paragraph{Characterisation of the family of steady-state solutions}
A biologically relevant steady-state solution of~\eqref{eq:KellerSegelSystem} subject to~\eqref{eq:ICs} and~\eqref{eq:BCs} is given by a pair of real, positive functions $u^\infty(x)$ and $c^\infty(x)$ that satisfy the following system of elliptic equations
\begin{equation}
\begin{cases}
\displaystyle{\nabla \cdot \left(\beta_u \ D(u^\infty) \ \nabla u^\infty - \chi \ \psi(u^\infty) \ u^\infty \ \nabla c^\infty \right) = 0, \quad u^\infty \equiv u^\infty(x),}
\\\\
\displaystyle{\beta_c \ \Delta c^\infty  + \alpha \ u^\infty - \kappa \ c^\infty = 0, \quad\quad\quad\quad\quad\quad\quad\quad\quad\;\;\; c^\infty \equiv c^\infty(x),}
\end{cases}
x \in \Omega \subset \mathbb{R}^d,
\label{eq:stationary-solutions-system}
\end{equation}
along with the boundary conditions
\begin{equation}\label{eq:BCsSS}
\nabla u^\infty(x) \cdot \hat{{\rm n}} = 0 \quad \text{and} \quad \nabla c^\infty(x) \cdot \hat{{\rm n}}, \quad x \in \partial \Omega,
\end{equation}
and the following integral identity, which follows from the principle of mass conservation,
\begin{equation} \label{eq:MCSS}
\int_{\Omega} u^\infty(x) \, \text{d}x = \int_{\Omega} u^0(x) \, \text{d}x =: M.
\end{equation}

The second equation in~\eqref{eq:stationary-solutions-system} along with the integral identity~\eqref{eq:MCSS} allow us to conclude that for a given value of $M$  there exists a unique positive homogeneous steady-state solution $\left(u^\infty, c^\infty \right) \equiv \left( u^\star, c^\star  \right)$, which is given by the pair
\begin{equation}
\label{eq:homSS}
\left( u^\star, c^\star  \right) = \left( \frac{M}{|\Omega|}, \frac{\alpha}{\kappa} \frac{M}{|\Omega|} \right).
\end{equation}

Moreover, since $D(u)$ and $\psi(u)$ are defined according to~\eqref{ass:volfil} and~\eqref{eq:exponential_chemofun}, the first equation in~\eqref{eq:stationary-solutions-system} along with the boundary conditions~\eqref{eq:BCsSS} give
\begin{equation}
\label{eq:stationary-solutions-equation}
\frac{D(u^\infty)}{u^\infty \ \psi(u^\infty)} \ \nabla u^\infty = \nu \ \nabla c^\infty \quad \Longrightarrow \quad u^\infty \ \exp\left( u^\infty / u_{\text{max}} \right) = \lambda \ \exp\left( \nu \ c^\infty\right),
\end{equation}
where $\nu := \chi/\beta_u$ and $\lambda$ is a real, positive constant that is uniquely defined by the integral identity~\eqref{eq:MCSS}. Building upon the analysis carried out in~\cite{Painter2002, Potapov2005}, we define $\Phi(z) := z \exp\left(z/u_\text{max}\right)$ and rewrite~\eqref{eq:stationary-solutions-equation} as
\begin{equation*}
\Phi(u^\infty) = \lambda \ \exp\left( \nu \ c^\infty\right).
\end{equation*}
Inverting the function $\Phi$, which is strictly increasing and thus invertible, we find
\begin{equation}
u^\infty = \Phi^{-1}\left( \lambda \exp\left( \nu \ c^\infty\right)\right) = u_\text{max} \ W \left(\frac{\lambda}{u_\text{max}} \exp\left( \nu \ c^\infty\right)\right),
\label{eq:uinfty}
\end{equation}
where $W$ is the multi-valued Lambert $W$ function~\cite{Corless1996}, which is implicitly defined by the relation
\begin{equation*}
W(z)  \exp\left(W(z) \right) = z, \quad \forall z \in \mathbb{C}.
\end{equation*}
Notice that $W(z)$ is real and single-valued for $z \in \mathbb{R}^+$ and, therefore, the right-hand side of~\eqref{eq:uinfty} is a real and single-valued function.
\begin{remark}
The Lambert $W$ function admits the following Taylor series expansion about the point $z=0$
\begin{equation*}
W(z) = z - z^2 + \frac{3}{2} z^3 + \mathcal{O}(z^4).
\end{equation*}
Hence, in the asymptotic regime $u_{\rm max} \to \infty$, i.e. when the generalised PKS model~\eqref{eq:KellerSegelSystem}-\eqref{eq:exponential_chemofun} formally reduces to the classical PKS system~\eqref{eq:KellerSegelSystemOri}, we can rewrite~\eqref{eq:uinfty} as
\begin{equation*}
W\left(\frac{\lambda}{u_{\rm max}}  \exp\left( \nu \ c^\infty\right)\right)  = \frac{\lambda}{u_{\rm max}} \exp\left( \nu \ c^\infty\right) - \frac{\lambda^2}{u_{\rm max}^2} \exp\left(2 \nu \ c^\infty\right) + h.o.t. \ .
\end{equation*}
Inserting the first order term of the above expansion into~\eqref{eq:uinfty} gives
\begin{equation*}
u^\infty = \lambda \exp\left( \nu \ c^\infty\right),
\end{equation*}
which is the well-known relation between $u^\infty$ and $c^\infty$ for the classical PKS model~\eqref{eq:KellerSegelSystemOri}.
\end{remark}

\begin{figure}[h!!]
		\centering
		\includegraphics[scale=0.09]{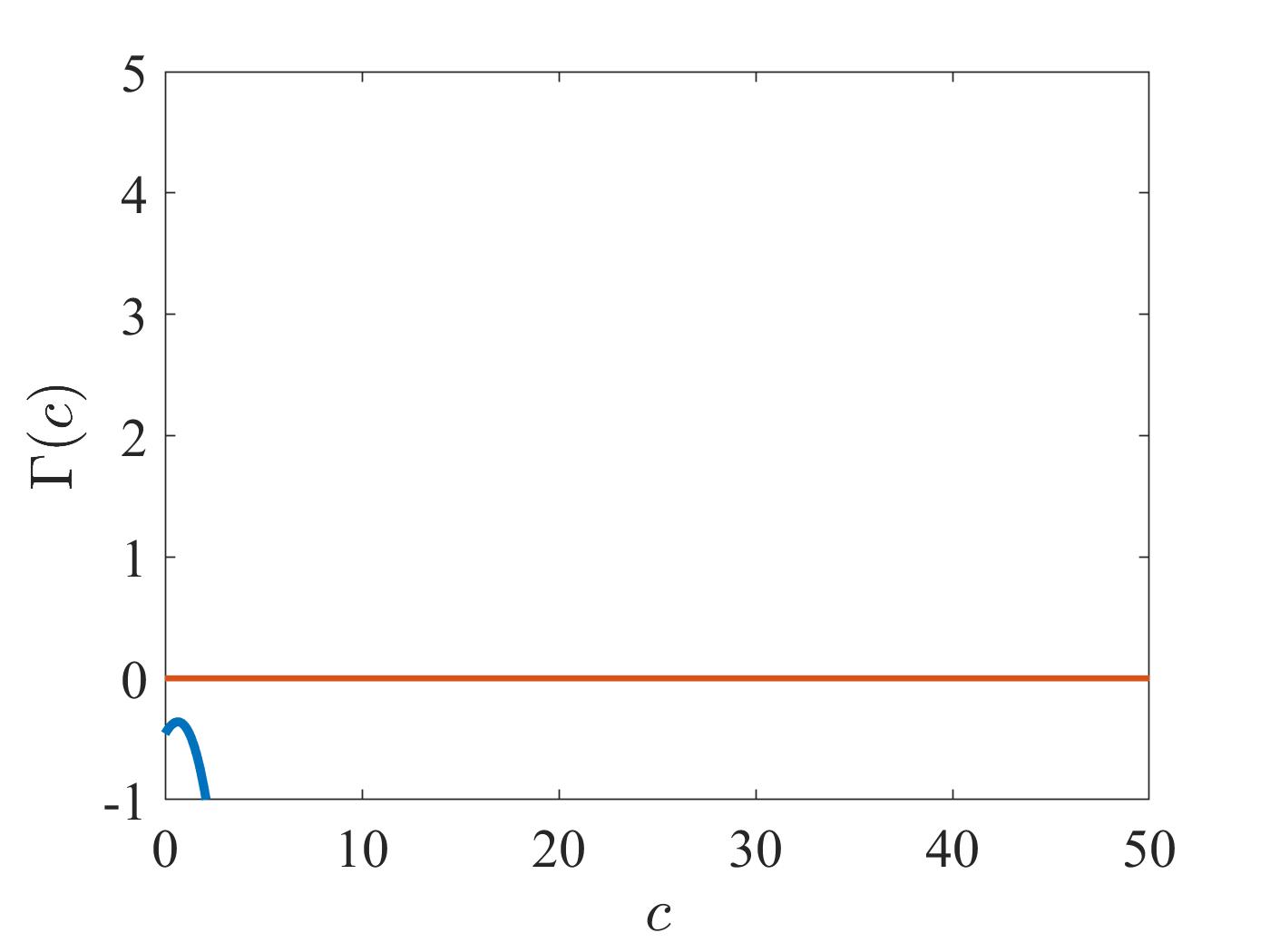}
		\includegraphics[scale=0.09]{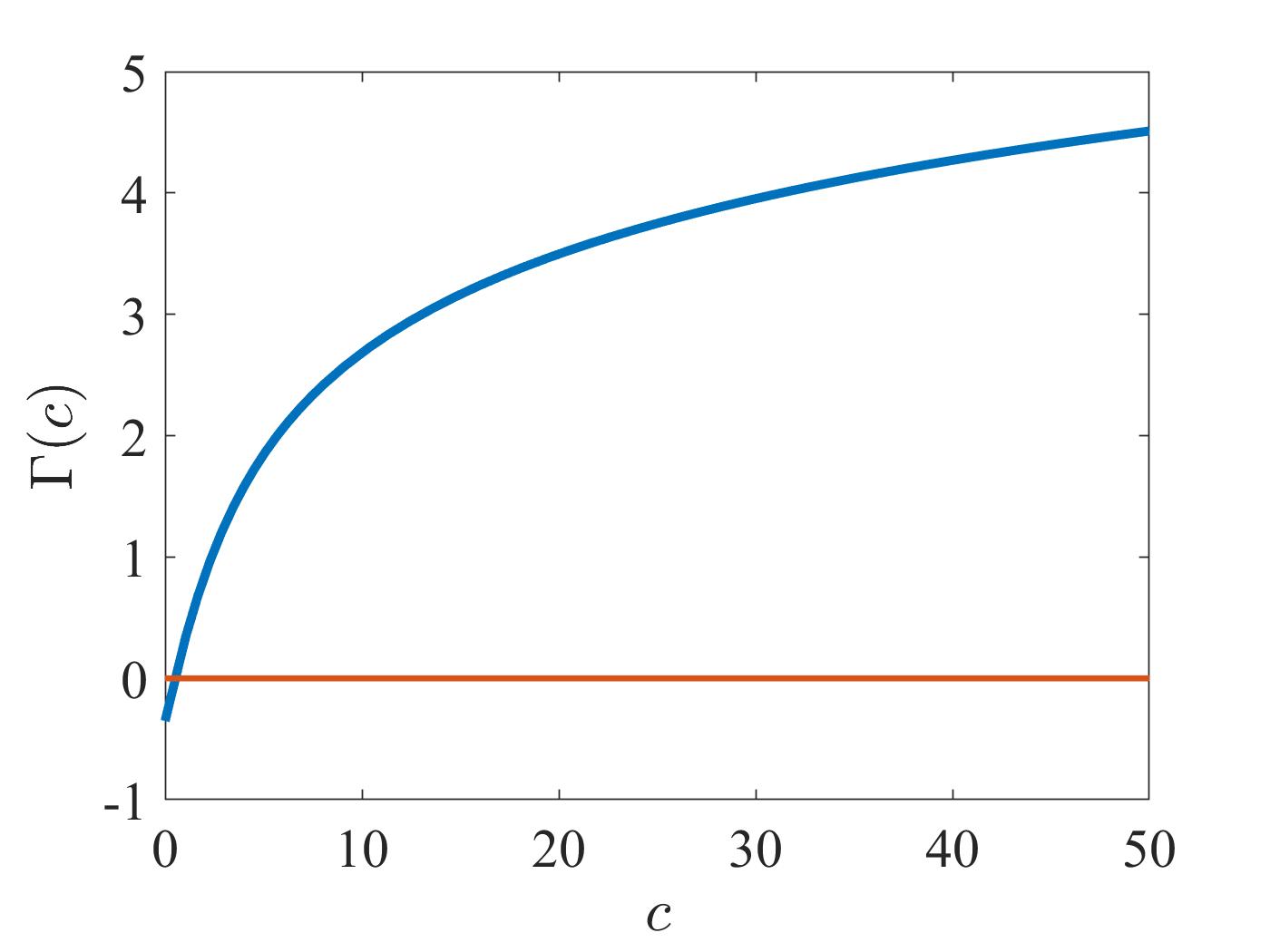}
		\includegraphics[scale=0.09]{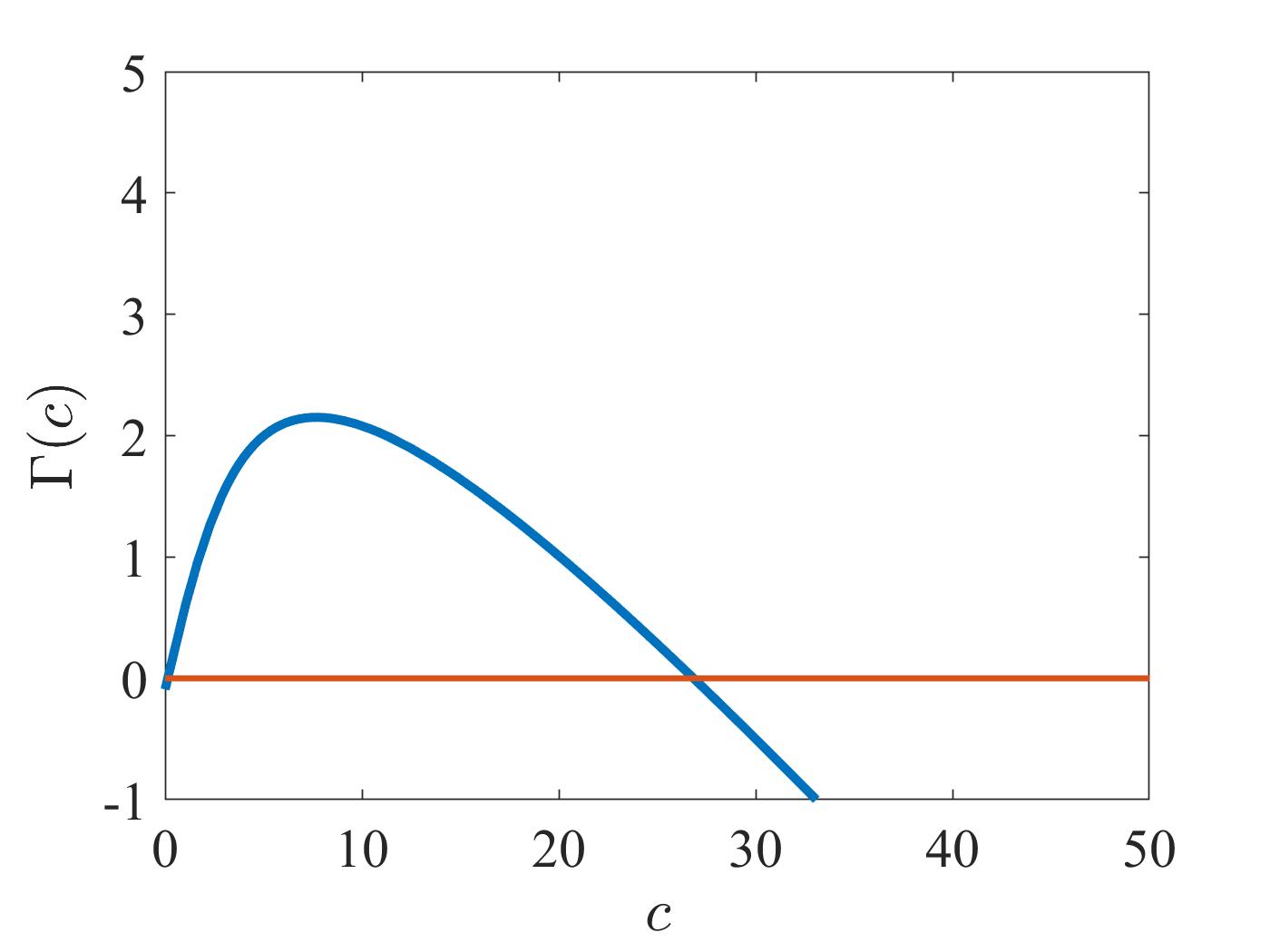}
		\caption{Sample plots of the function $\Gamma(c)$ for three different sets of values of the parameters in~\eqref{defgamma}.}
		\label{fig:roots}
	\end{figure}

When $\Omega := (0,L) \subset \mathbb{R}$ with $L>0$, substituting~\eqref{eq:uinfty} into the second equation in~\eqref{eq:stationary-solutions-system}, introducing the notation
\begin{equation}
\label{defgamma}
\Gamma(c) := \kappa \, c - \alpha \, u_\text{max} \, W \left(\frac{\lambda}{u_\text{max}} \exp\left( \nu \ c\right)\right)
\end{equation}
and imposing the boundary conditions~\eqref{eq:BCsSS} gives the following second order differential equation 
\begin{equation*}
\beta_c \ \dfrac{{\rm d}^2 c^{\infty}}{{\rm d} x^2} = \Gamma(c^{\infty}), \quad c^\infty \equiv c^\infty(x), \quad x \in (0,L),
\end{equation*}
subject to
\begin{equation*}
\dfrac{{\rm d} c^{\infty}(0)}{{\rm d} x} = 0 \quad \text{and} \quad \dfrac{{\rm d} c^{\infty}(L)}{{\rm d} x} = 0.
\end{equation*}
As similarly done in~\cite{Painter2002}, further insight into the properties of $c^{\infty}(x)$ can be gained using phase-plane methods. We rewrite the latter second order differential equation as the following system of first order differential equations
\begin{equation}
\label{eq:cinfty_eqODEs}
\begin{cases}
\dfrac{{\rm d} c^{\infty}}{{\rm d} x} = w, \quad\quad\quad\,\, c^\infty \equiv c^\infty(x),\\\\
\dfrac{{\rm d} w}{{\rm d} x} = \dfrac{1}{\beta_c}\Gamma(c^{\infty}), \quad w \equiv w(x),
\end{cases}
x \in (0,L).
\end{equation}
In contrast to the case considered in~\cite{Painter2002}, here we cannot determine the number of equilibria of~\eqref{eq:cinfty_eqODEs} due to the fact that we cannot infer the number of roots of $\Gamma$. In particular, the condition $\Gamma(0) < 0$ can be deduced from the fact that the function $W$ is increasing, but we do not have enough information to characterise the behaviour of $\Gamma(c^{\infty})$ as $c^{\infty} \to \infty$. However, as summarised by the plots in Figure~\ref{fig:roots}, numerical simulations indicate that, depending on the values of the parameters in~\eqref{defgamma}, we can generally expect zero, one or two distinct non-negative roots. 

When $c_1 \in \mathbb{R}^+$ is the unique non-negative root of the function $\Gamma$, we have that $\Gamma^{\prime}(c_1) > 0$ and, therefore, the corresponding equilibrium $(c_1,0)$ of~\eqref{eq:cinfty_eqODEs} is a saddle point and the positive homogeneous steady-state state is stable. On the other hand, if the function $\Gamma$ has two non-negative roots $c_1 \in \mathbb{R}^+$ and $c_2 \in \mathbb{R}^+$ with $c_2>c_1$, then necessarily $\Gamma^{\prime}(c_1) > 0$ and $\Gamma^{\prime}(c_2) < 0$. Hence, the corresponding equilibria $(c_1,0)$ and $(c_2,0)$ of~\eqref{eq:cinfty_eqODEs} will be a saddle point and a centre, respectively, and in this case we can expect steady state solutions with multiple peaks.

\paragraph*{Linear stability analysis of the positive homogeneous steady state}
It is straightforward to show that the steady state $\left( u^\star, c^\star  \right)$ is linearly stable to spatially homogeneous perturbations. Furthermore, in order to study the linear stability of the steady state $\left( u^\star, c^\star  \right)$ to spatially inhomogeneous perturbations, we make the ansatz
\begin{equation}
u(t,x) = u^\star + \tilde{u} \ \exp{(\sigma t)} \ \varphi_k(x), \quad c(t,x) = c^\star + \ \tilde{c} \ \exp{(\sigma t)} \ \varphi_k(x),
\label{eq:linearised_solutions}
\end{equation}
where $\tilde{u},\tilde{c} \in \mathbb{R}_*$ with $|\tilde{u}| \ll 1$ and $|\tilde{c}| \ll 1$, $\sigma \in \mathbb{C}$ and $\{\varphi_k\}_{k \geq 1}$ are the eigenfunctions of the Laplace operator indexed by the wavenumber $k$. 
Linearising~\eqref{eq:KellerSegelSystem} about the steady state $(u^\star,c^\star)$ and substituting~\eqref{eq:linearised_solutions} into the linearised system of equations yields  
\begin{equation*}
\begin{cases}
\sigma \ \tilde{u} = -k^2 \ \beta_u \ D(u^\star) \ \tilde{u} + k^2 \ \chi \ \psi(u^\star) \ u^\star \ \tilde{c}, \\
\sigma \ \tilde{c} = -k^2  \  \beta_c  \ \tilde{c} + \alpha  \  \tilde{u} - \kappa  \  \tilde{c}.
\end{cases}
\end{equation*}
For the above system to admit a solution $(\tilde{u},\tilde{c}) \in \mathbb{R}^2_*$ we need
\begin{equation*}
\sigma^2 + \left[k^2 (\beta_u D(u^\star) + \beta_c) + \kappa \right] \sigma + k^4 \beta_u \beta_c D(u^\star) + k^2 \left(\kappa \beta_u D(u^\star)-\alpha \chi u^\star \psi(u^\star) \right) = 0.
\end{equation*}
The steady state $\left( u^\star, c^\star  \right)$ will be driven unstable by spatially inhomogeneous perturbations (\emph{i.e.} spatial patterns will emerge) if ${\rm Re}(\sigma)>0$. Since $k^2 (\beta_u D(u^\star) + \beta_c) + \kappa>0$, from the above polynomial equation for $\sigma$ we conclude that ${\rm Re}(\sigma)>0$ if 
\begin{equation*}
k^4 \beta_u \beta_c D(u^\star) + k^2 \left(\kappa \beta_u D(u^\star)-\alpha \chi u^\star \psi(u^\star) \right) < 0
\end{equation*}
for an interval of values of $k^2 \in \mathbb{R}^+_*$, that is, if
\begin{equation*}
\alpha \chi u^\star \psi(u^\star) -\kappa \beta_u D(u^\star) > 0 \quad \Longrightarrow \quad \chi > \frac{\kappa \beta_u D(u^\star)}{\alpha u^\star \psi(u^\star)}
\end{equation*}
and
\begin{equation*}
0 < k^2 < \frac{\alpha \chi u^\star \psi(u^\star)-\kappa \beta_u D(u^\star)}{\beta_u \beta_c D(u^\star)}.
\end{equation*}
In the case where the functions $D(u)$ and $\psi(u)$ are defined according to~\eqref{ass:volfil} and~\eqref{eq:exponential_chemofun}, respectively, the above conditions reduce to
\begin{equation}
\label{eq:condinst}
\chi > \frac{\kappa \beta_u (1+u^\star/u_{\text{max}})}{\alpha u^\star}
\end{equation}
and
\begin{equation}
0 < k^2 <  k^2_\text{max} \quad \text{with} \quad k^2_\text{max} := \frac{\alpha \chi u^\star - \kappa \beta_u \left(1 + u^\star/u_{\text{max}} \right)}{\beta_u \beta_c \left(1 + u^\star /u_{\text{max}}\right)}.
\label{eq:modes}
\end{equation}
Notice that $k^2_\text{max}$ is an increasing function of the chemotactic sensitivity $\chi$, which implies that if $\chi$ increases then the most unstable mode associated with the largest eigenvalue $\sigma$ increases and the range of unstable modes broadens~\cite{Potapov2005}. This is confirmed by the numerical results presented in Section~\ref{sec:comparison}\ref{sec:numres}.

\subsection{Main results of numerical simulations}
\label{sec:numres}
First, we present the results of  base-case numerical simulations showing sample spatial patterns of the discrete model and the generalised PKS model~\eqref{eq:KellerSegelSystem}-\eqref{eq:exponential_chemofun} (\emph{i.e.} the system of coupled parabolic equations~\eqref{modexpli}). Then, we investigate how the spatial patterns produced by the two models can vary with the strength of chemotactic sensitivity (\emph{i.e.} the value of the parameter $\eta$ of the discrete model and the value of the corresponding parameter $\chi$ of the continuum model), the size of the cell population, \emph{i.e.} the quantities
\begin{equation*}
\sum_i n_i^0 \quad \text{and} \quad \int_{\Omega} u^0(x) \, {\rm d}x,
\end{equation*}
and the critical cell density $u_{\rm max}$ in definition~\eqref{eq:exponential_chemofun}. Finally, we explore the existence of scenarios in which differences between spatial patterns produced by the two models can emerge due to effects associated with small cell numbers, which reduce the quality of the approximation of the discrete model provided by the continuum model. A complete description of the set-up of numerical simulations, the algorithmic rules that underlie computational simulations of the discrete model, and the numerical methods used to solve numerically the generalised PKS model~\eqref{eq:KellerSegelSystem}-\eqref{eq:exponential_chemofun} are given in Appendices~B and~. In particular, for all the numerical simulations we report on in this section, the parameter values are such that either conditions~\eqref{ass:quottozero1d} or conditions~\eqref{ass:quottozero2d} are satisfied.

\begin{figure}[!h]
	\centering
	\includegraphics[width=\textwidth]{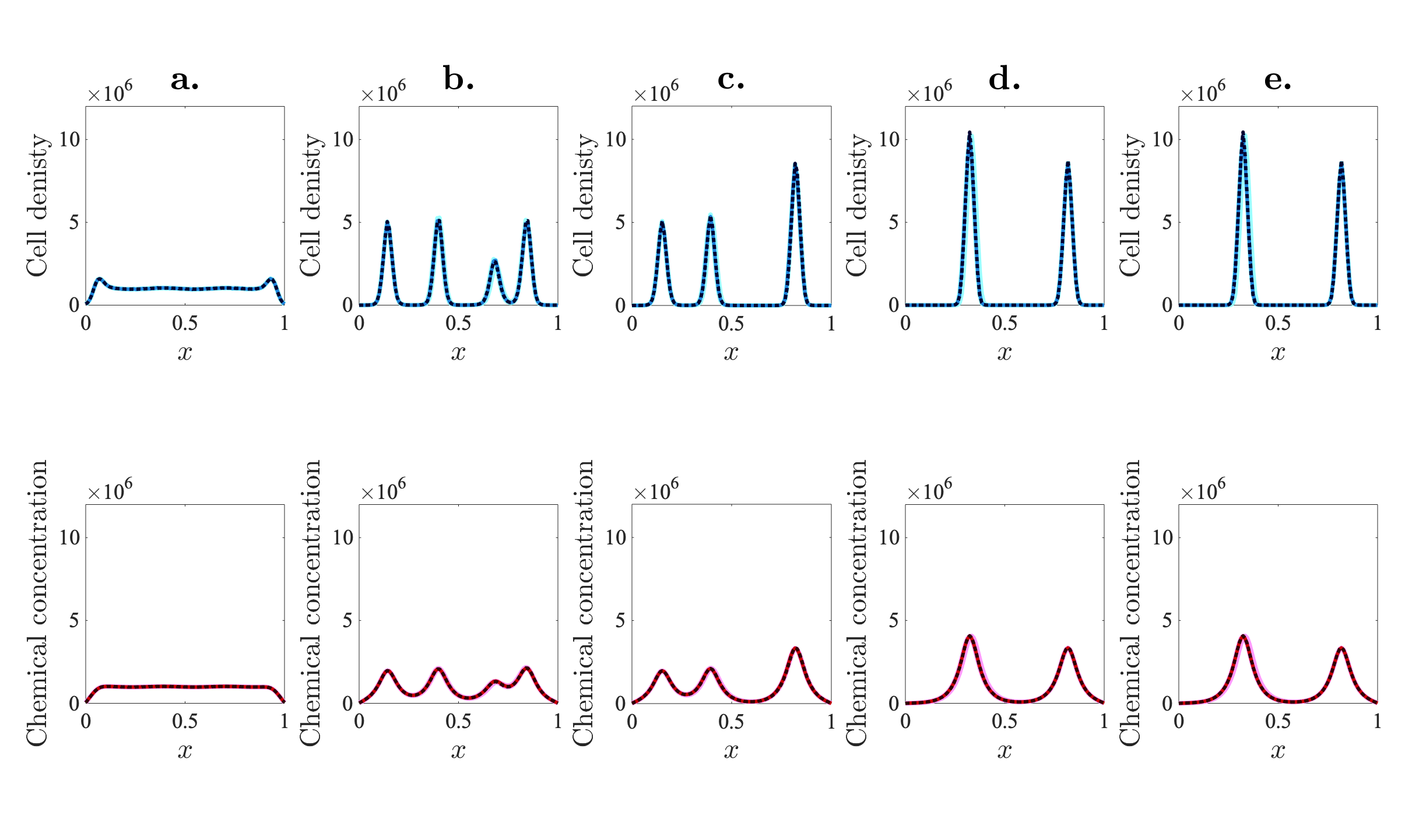}\\
	\includegraphics[width=0.7\textwidth]{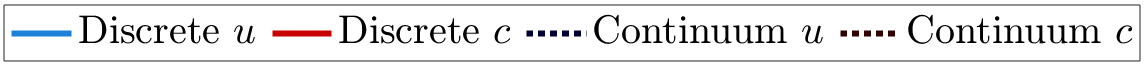}
\caption{ {\bf Base-case numerical results in one dimension} Comparison between the simulation results for the discrete model (solid lines) and numerical solutions of the generalised PKS model~\eqref{eq:KellerSegelSystem}-\eqref{eq:exponential_chemofun} (dotted lines). The top panels display the cell density (blue lines) and the bottom panels display the concentration of chemoattractant (red lines) at five successive time instants -- \emph{i.e.} {\bf a.} $t=$ 1, {\bf b.} $t=$ 25, {\bf c.} $t=$ 50, {\bf d.} $t=$ 300, {\bf e.} $t=$ 500. The results from the discrete model correspond to the average over five realisations of the underlying biased random walk. The cell density and the concentration of chemoattractant resulting from each realisation are plotted in pale blue and magenta, respectively, to demonstrate the robustness of the results obtained. A complete description of the set-up of numerical simulations and the numerical methods employed is given in Appendix~B.}
	\label{fig:basecase}
\end{figure}

\begin{figure}[!h]
	\centering
	\includegraphics[width=\textwidth]{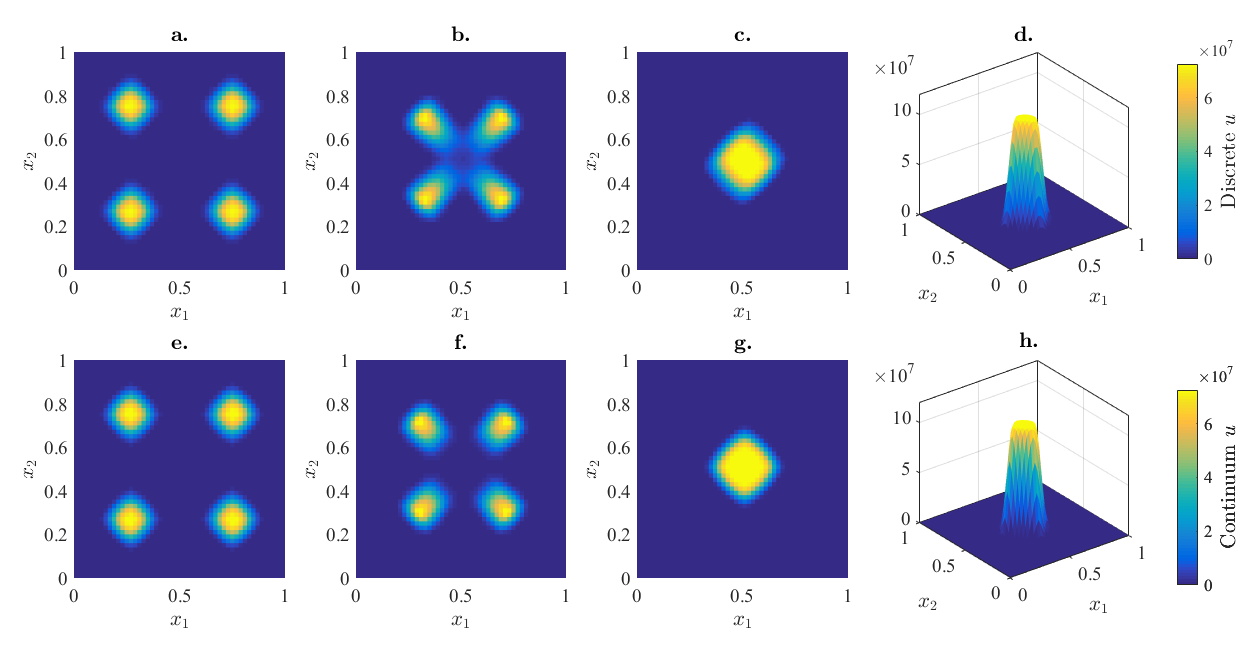}
\caption{ {\bf Base-case numerical results in two dimensions} Comparison between the simulation results for the discrete model and numerical solutions of the generalised PKS model~\eqref{eq:KellerSegelSystem}-\eqref{eq:exponential_chemofun}. Panels {\bf a.-c.} display the discrete cell density at three successive time instants -- \emph{i.e.} {\bf a.} $t=$ 0, {\bf b.} $t=$ 5, {\bf c.} $t=$ 15 -- while panels {\bf e.-g.} display the corresponding solutions of the generalised PKS model. Panels {\bf d., h.} display a side-on view of the cell density at the end of numerical simulations for the discrete model and the generalised PKS model, respectively. The results from the discrete model correspond to the average over two realisations of the underlying biased random walk. A complete description of the set-up of numerical simulations and the numerical methods employed is given in Appendix~C.}
	\label{fig:basecase2D}
\end{figure}

\paragraph{Base-case numerical results} Figure~\ref{fig:basecase} and Figure~\ref{fig:basecase2D} demonstrate that there is an excellent quantitative match between numerical solutions of the generalised PKS model~\eqref{eq:KellerSegelSystem}-\eqref{eq:exponential_chemofun} and the results of numerical simulations of the discrete model, both in one and in two spatial dimensions. In agreement with the results of linear stability analysis carried out in Section~\ref{sec:comparison}\ref{sec:stability}, since condition~\eqref{eq:condinst} is satisfied under the parameter setting considered here, spatial patterns are formed. In the one-dimensional case ({\it cf.} Figure~\ref{fig:basecase}), we first observe the emergence of four peaks in the cell density, as well as in the concentration of chemoattractant, which then merge into three peaks before coalescing into two peaks. On the other hand, in the two-dimensional case ({\it cf.} Figure~\ref{fig:basecase2D}), we observe the emergence of a plateau. This is due to the interplay between the tendency of cell density to become locally concentrated and the fact that the type of nonlinear diffusion considered ensures boundedness of the solutions. Later in this section, we further investigate how stationary solutions are affected by the size of the cell population. 

\begin{figure}[!h]
	\centering
	\includegraphics[width=\textwidth]{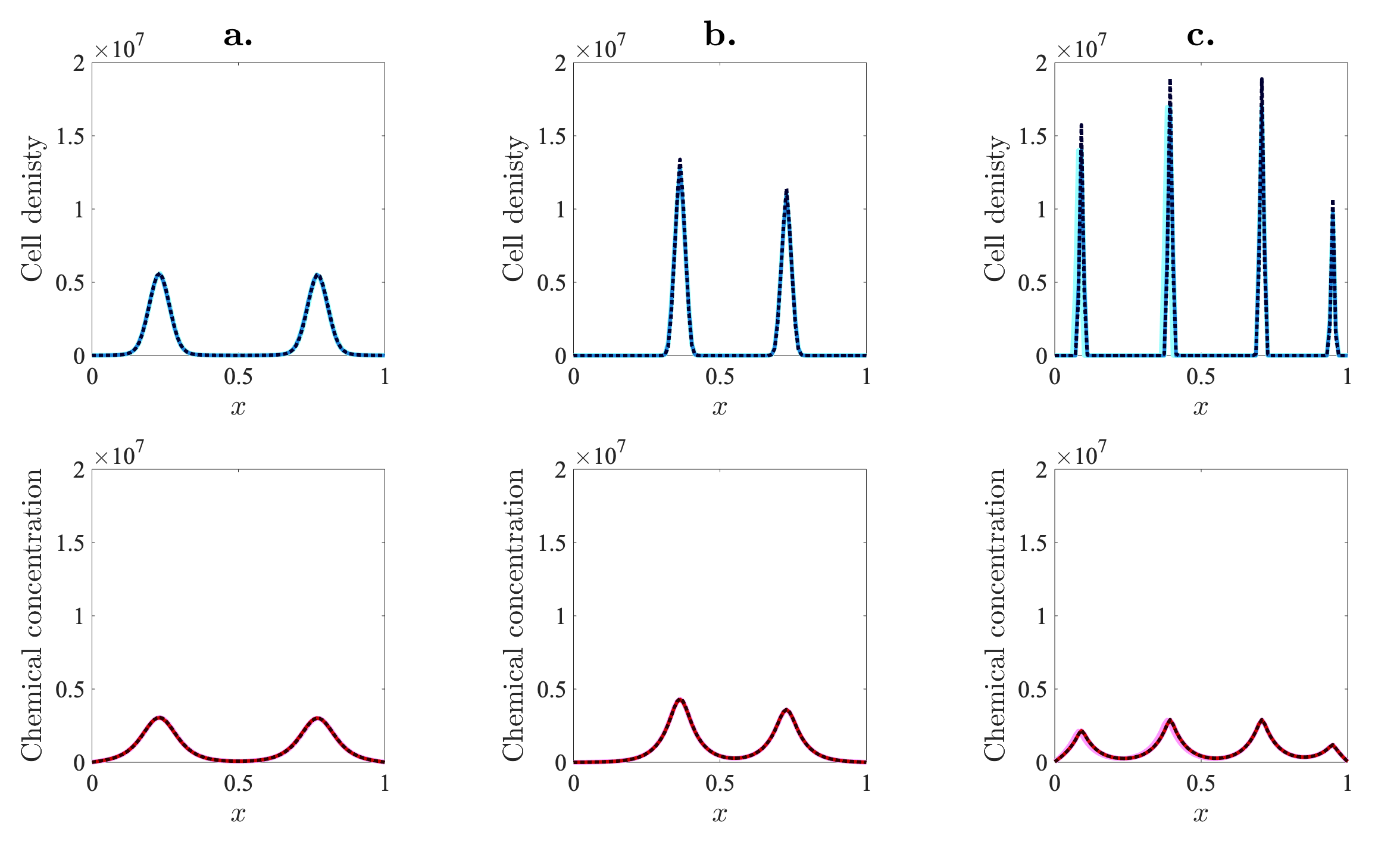}\\
	\includegraphics[width=0.7\textwidth]{Legend1.png}
\caption{ {\bf Effect of the strength of chemotactic sensitivity} Comparison between the simulation results for the discrete model (solid lines) and numerical solutions of the generalised PKS model~\eqref{eq:KellerSegelSystem}-\eqref{eq:exponential_chemofun} (dotted lines). The top panels display the cell density (blue lines) and the bottom panels display the concentration of chemoattractant (red lines) at the end of numerical simulations (\emph{i.e.} at numerical equilibrium). The three sets of panels refer to different values of the parameter $\eta$ in the discrete model -- \emph{i.e.} {\bf a.} $\eta=$ 0.9801, {\bf b.} $\eta=$ 4.9005, {\bf c.} $\eta=$ 294.03 -- which correspond to different values of the parameter $\chi$ defined via~\eqref{ass:quottozero1d} in the continuum model. The results from the discrete model correspond to the average over five realisations of the underlying biased random walk. The cell density and the concentration of chemoattractant resulting from each realisation are plotted in pale blue and magenta, respectively, to demonstrate the robustness of the results obtained. A complete description of the set-up of numerical simulations and the numerical methods employed is given in Appendix~B.}
	\label{fig:chi}
\end{figure}

\paragraph{Effect of the strength of chemotactic sensitivity} 
Figure~\ref{fig:chi} indicates that, coherently with relation~\eqref{eq:modes}, the number of peaks observed at numerical equilibrium increases with the value of the parameter $\eta$ in the discrete model and the corresponding value of the parameter $\chi$ defined via~\eqref{ass:quottozero1d} in the continuum model. For all values of $\eta$ considered, the numerical results obtained indicate excellent agreement between the simulation results for the discrete model and numerical solutions of the generalised PKS model~\eqref{eq:KellerSegelSystem}-\eqref{eq:exponential_chemofun}.

\begin{figure}[!h]
\centering
	\includegraphics[width=\textwidth]{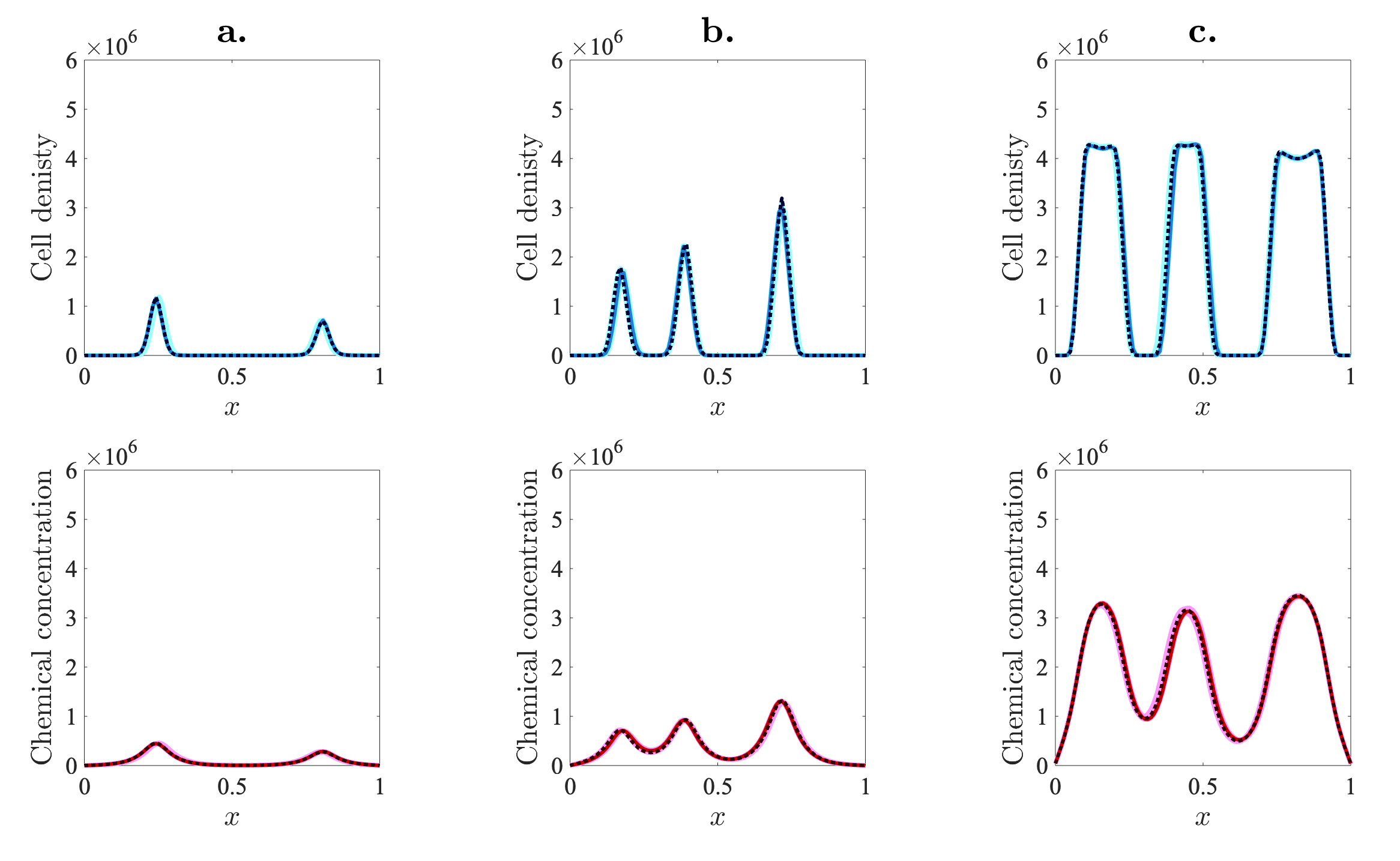}\\
	\includegraphics[width=0.7\textwidth]{Legend1.png}
\caption{ {\bf Effect of the size of the cell population in one dimension} Comparison between the simulation results for the discrete model (solid lines) and numerical solutions of the generalised PKS model~\eqref{eq:KellerSegelSystem}-\eqref{eq:exponential_chemofun} (dotted lines). The top panels display the cell density (blue lines) and the bottom panels display the concentration of chemoattractant (red lines) at the end of numerical simulations for different values of the size of the cell population, that is, different values of the parameter $B$ in~\eqref{eq:ICB} -- \emph{i.e.} {\bf a.} $B=$ 0.25, {\bf b.} $B=$ 1, {\bf c.} $B=$ 5. The results from the discrete model correspond to the average over five realisations of the underlying biased random walk. The cell density and the concentration of chemoattractant resulting from each realisation are plotted in pale blue and magenta, respectively, to demonstrate the robustness of the results obtained. A complete description of the set-up of numerical simulations and the numerical methods employed is given in Appendix~B.}
	\label{fig:initial}
\end{figure}

\begin{figure}[!h]
\centering
	\includegraphics[width=\textwidth]{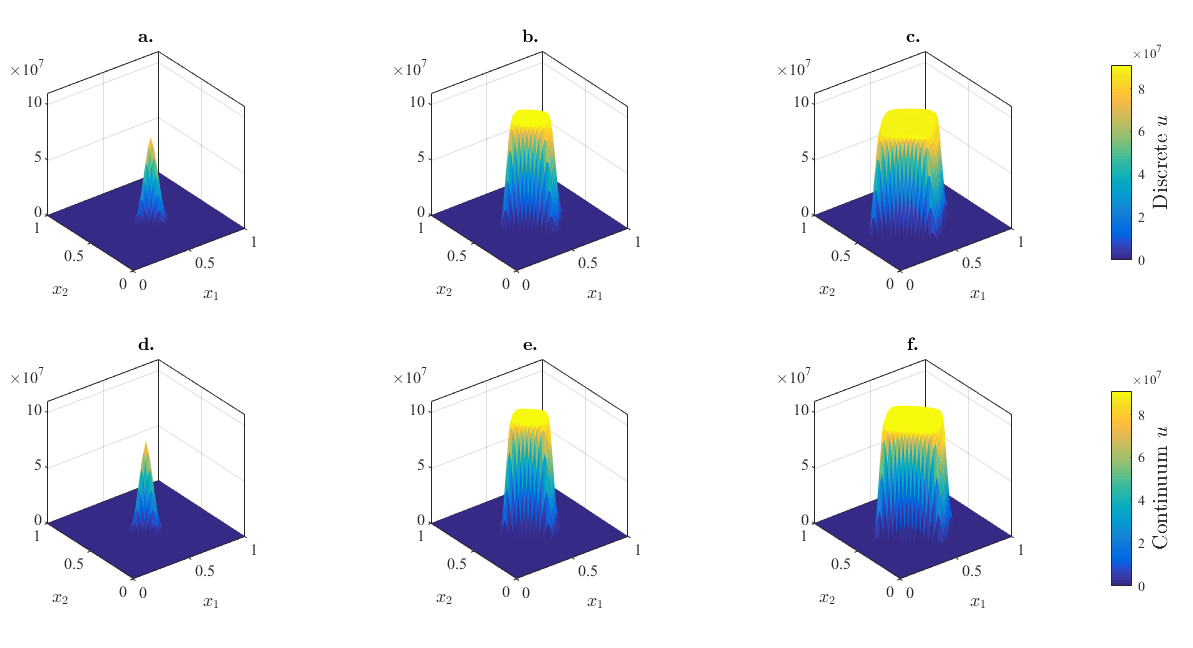}
\caption{ {\bf Effect of the size of the cell population in two dimensions} Comparison between the simulation results for the discrete model ({\bf a.-c.}) and numerical solutions of the generalised PKS model~\eqref{eq:KellerSegelSystem}-\eqref{eq:exponential_chemofun} ({\bf d.-f.}). The different panels display the cell density at the end of numerical simulations for different values of the size of the cell population, that is, different values of the parameter $B$ in~\eqref{eq:ICB} -- \emph{i.e.} {\bf a.-d.} $B=$ 0.1, {\bf b.-e.} $B=$ 1, {\bf c.-f.} $B=$ 2. The results from the discrete model in panels {\bf a.-b.} correspond to the average over two realisations of the underlying biased random walk, while the results in panel {\bf c.} correspond to a single realisation. A complete description of the set-up of numerical simulations and the numerical methods employed is given in Appendix~C.}
	\label{fig:initial2D}
\end{figure}

\paragraph{Effect of the size of the cell population}
Figures~\ref{fig:initial} and~\ref{fig:initial2D} display the plots of the cell density and the concentration of chemoattractant obtained at the end of numerical simulations for different values of the size of the cell population, \emph{i.e.} considering initial conditions $n_i^0$ and $u^0(x)$ such that
\begin{equation}
\label{eq:ICB}
\sum_i n_i^0 = B \, M  \quad \text{and} \quad \int_{\Omega} u^0(x) \, {\rm d}x = B \, M,
\end{equation}
with $M>0$ fixed and for different values of $B > 0$. These plots in Figure~\ref{fig:initial} show that, in contrast to the classical PKS model, incorporating volume-filling effects through definitions~\eqref{ass:volfil} and \eqref{eq:exponential_chemofun} prevents unphysical finite-time blow-up. Moreover, the spatial patterns produced by the two models change with the size of the cell population and, for all values of $B$ in~\eqref{eq:ICB} considered, there is an excellent quantitative match between the results for numerical simulations of the discrete model and numerical solutions of the generalised PKS model~\eqref{eq:KellerSegelSystem}-\eqref{eq:exponential_chemofun}, both in one and in two spatial dimensions.

\begin{figure}[!h]
	\centering
	\includegraphics[width=\textwidth]{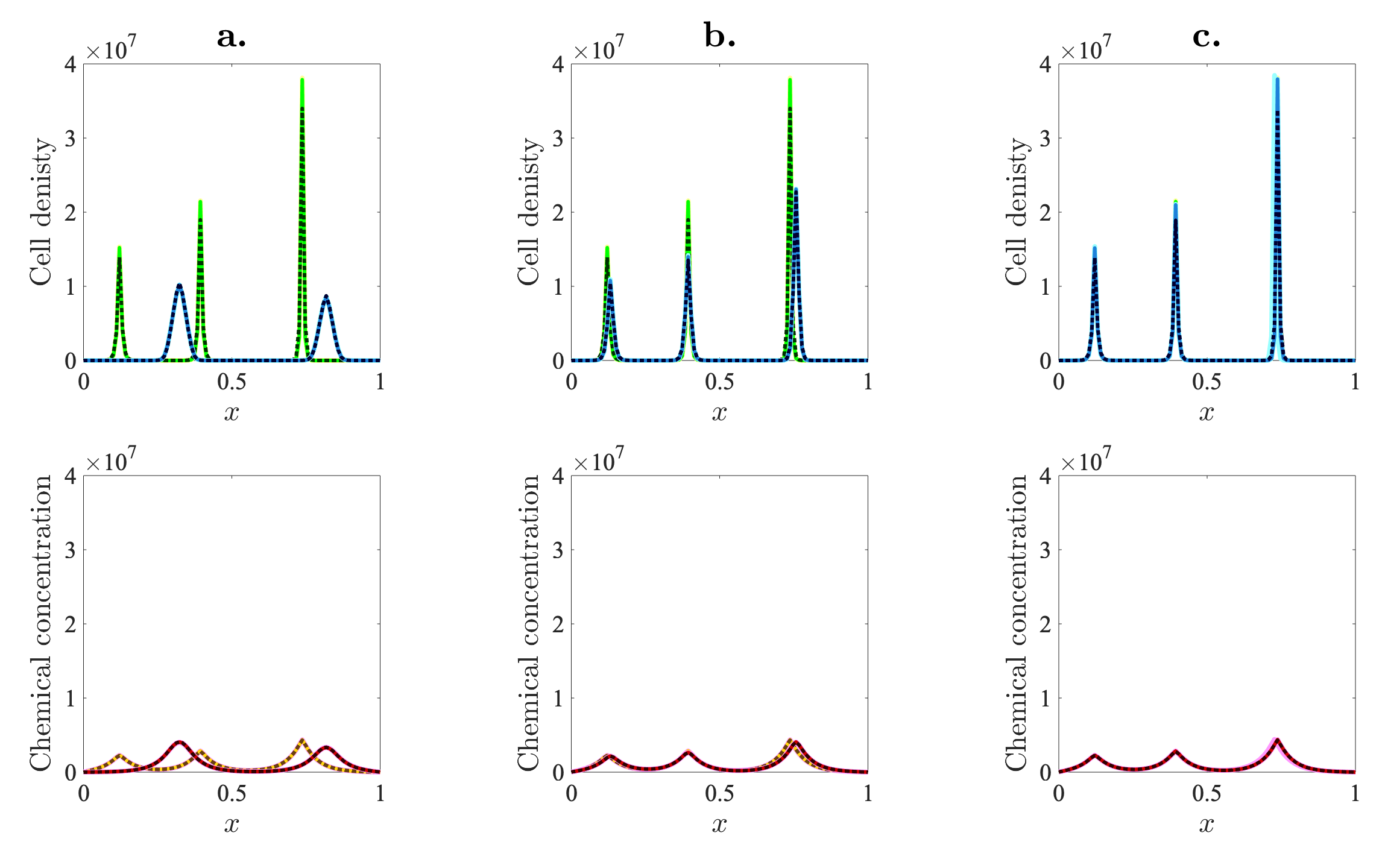}\\
	\includegraphics[width=\textwidth]{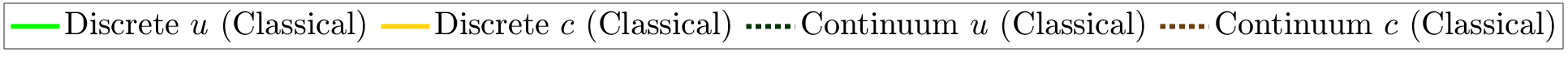}\\
	\includegraphics[width=\textwidth]{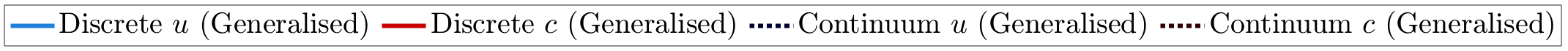}\\
\caption{ {\bf Effect of the critical cell density $u_{\rm max}$ in one dimension} Comparison between numerical solutions of the generalised PKS model~\eqref{eq:KellerSegelSystem}-\eqref{eq:exponential_chemofun}, numerical solutions of the classical PKS model~\eqref{eq:KellerSegelSystemOri}, and the simulation results for the corresponding discrete models -- \emph{i.e.} the discrete model with the function $\psi$ in~\eqref{e:left}, \eqref{e:right}, \eqref{e:leftdif} and \eqref{e:rightdiff} defined via~\eqref{eq:exponential_chemofun}  or with $\psi \equiv 1$, respectively. The solid, blue lines and the solid, red lines highlight the cell density and the concentration of chemoattractant at the end of numerical simulations of the discrete model with the function $\psi$ defined via~\eqref{eq:exponential_chemofun} (\emph{i.e.} with volume-filling effects). On the other hand, the solid, green lines and the solid, yellow lines highlight the cell density and the concentration of chemoattractant at the end of numerical simulations of the discrete model with $\psi \equiv 1$ (\emph{i.e.} without volume-filling effects). The dotted lines highlight the numerical solutions of the corresponding PKS models. Different panels refer to different values of the parameter $u_{\rm max}$ -- \emph{i.e.} {\bf a.} $u_{\rm max}=2\times 10^{6}$, {\bf b.} $u_{\rm max}=2\times 10^{7}$, {\bf c.} $u_{\rm max}=2\times 10^{9}$. The results from the discrete model with the function $\psi$ defined via~\eqref{eq:exponential_chemofun} correspond to the average over thirty realisations of the underlying biased random walk, while the results from the discrete model with $\psi \equiv 1$ correspond to the average over ten realisations. The cell density and the concentration of chemoattractant resulting from each realisation are plotted in paler colours to demonstrate the robustness of the results obtained. A complete description of the set-up of numerical simulations and the numerical methods employed is given in Appendix~B.}
	\label{fig:infinity}
\end{figure}

\begin{figure}[!h]
	\centering
	\includegraphics[width=\textwidth]{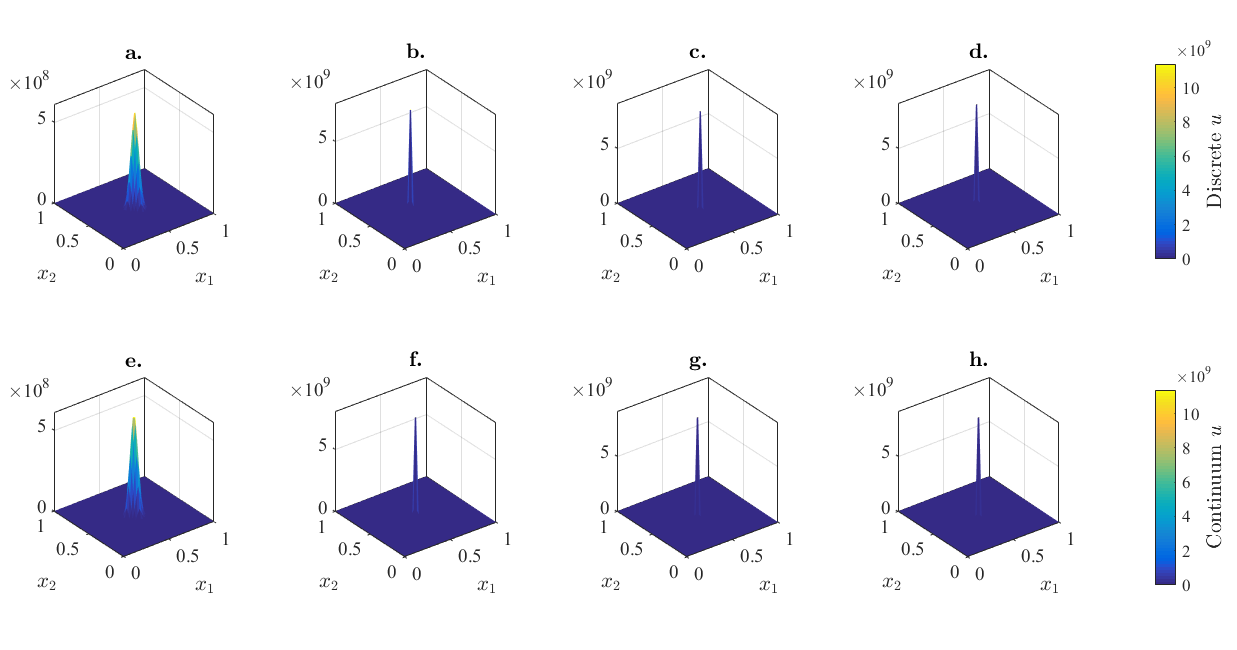}
\caption{ {\bf Effect of the critical cell density $u_{\rm max}$ in two dimensions} Comparison between numerical solutions of the generalised PKS model~\eqref{eq:KellerSegelSystem}-\eqref{eq:exponential_chemofun}, numerical solutions of the classical PKS model~\eqref{eq:KellerSegelSystemOri}, and the simulation results for the corresponding discrete models -- \emph{i.e.} the discrete model with the function $\psi$ in the probabilities of cell movement defined via~\eqref{eq:exponential_chemofun} or with $\psi \equiv 1$, respectively. Panels {\bf a.-c.} display the cell density at the end of numerical simulations of the discrete model with the function $\psi$ defined via~\eqref{eq:exponential_chemofun} (\emph{i.e.} with volume-filling effects) for different values of the parameter $u_{\rm max}$ - \emph{i.e.} {\bf a.} $u_{\rm max}=1\times 10^{8}$, {\bf b.} $u_{\rm max}=1\times 10^{10}$, {\bf c.} $u_{\rm max}=1\times 10^{11}$. The corresponding numerical solutions of the generalised PKS model~\eqref{eq:KellerSegelSystem}-\eqref{eq:exponential_chemofun} are displayed in panels {\bf e.-g.}. Panel {\bf d.} displays the cell density at the end of numerical simulations of the discrete model with $\psi \equiv 1$ (\emph{i.e.} without volume-filling effects), and the corresponding numerical solution of the classical PKS model~\eqref{eq:KellerSegelSystemOri} is displayed in panel {\bf h.}. The results displayed in panel {\bf a.} correspond to the average over two realisations of the underlying biased random walk, while all other results from the discrete models correspond to a single realisation. A complete description of the set-up of numerical simulations and the numerical methods employed is given in Appendix~C.}	
	\label{fig:infinity2D}
\end{figure}

\paragraph{Effect of the critical cell density $u_{\rm max}$}
The generalised PKS model~\eqref{eq:KellerSegelSystem}-\eqref{eq:exponential_chemofun} formally reduces to the classical PKS model~\eqref{eq:KellerSegelSystemOri} as $u_{\text{max}} \to \infty$, since $\psi(u)\equiv1$ when $u_{\text{max}} \to \infty$ and, therefore, $D(u)\equiv1$ and $\psi(u)\equiv1$ in this limit. As a result, we expect the generalised PKS model, and thus our discrete model, to exhibit the same spatial patterns as those produced by the classical PKS model in such an asymptotic regime. This is confirmed by the numerical results presented in Figures~\ref{fig:infinity} and Figure~\ref{fig:infinity2D}. These results demonstrate how increasing values of $u_{\text{max}}$ lead to a better match between numerical solutions of the generalised PKS model and those of the classical PKS model, both in one and in two spatial dimensions (\emph{cf.} Figures~\ref{fig:infinity}{\bf a., b.} and Figures~\ref{fig:infinity2D}{\bf a.-c., e.-g.}), and a perfect quantitative match is ultimately obtained for $u_{\text{max}}$ sufficiently high (\emph{cf.} Figures~\ref{fig:infinity}{\bf c.} and~\ref{fig:infinity2D}{\bf h.}). In all cases, there is an excellent agreement between numerical solutions of the PKS models and the results for numerical simulations of the corresponding discrete models. Notice that the discrete model corresponding to the classical PKS model~\eqref{eq:KellerSegelSystemOri} is defined by assuming $\psi \equiv 1$ in~\eqref{e:left}, \eqref{e:right}, \eqref{e:leftdif} and \eqref{e:rightdiff}, and in their two-dimensional counterparts (\emph{i.e.} there are no volume-filling effects). We expect analogous results to hold in higher spatial dimensions (\emph{i.e.} when $d \geq 3$) whereby the solutions of the classical PKS model are known to blow up for cell populations of arbitrarily small size~\cite{Winkler2013}. 

\begin{figure}[!h]
	\centering
	\includegraphics[width=1\textwidth]{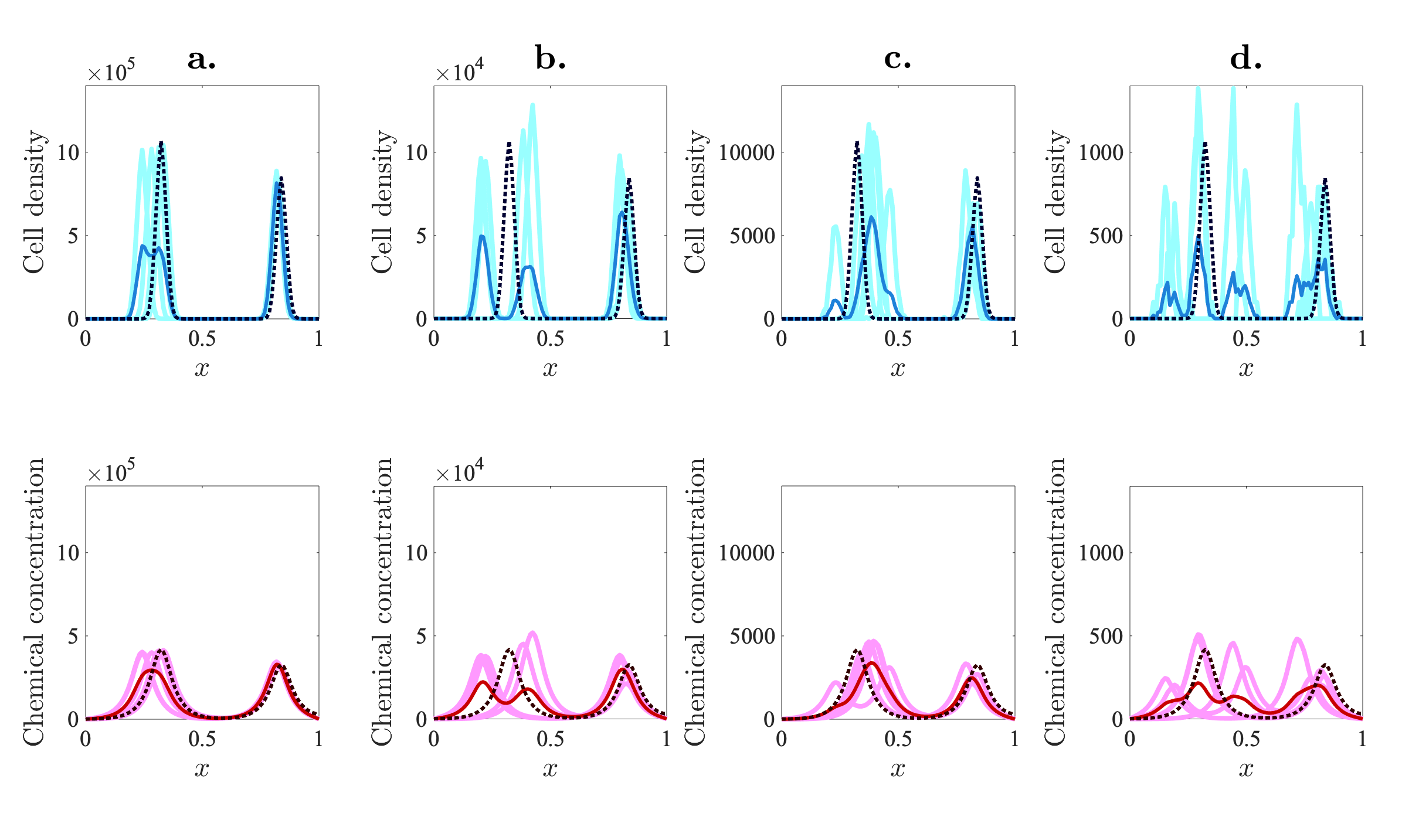}\\
	\includegraphics[width=0.7\textwidth]{Legend1.png}
\caption{{\bf Emergence of differences between spatial patterns produced by the discrete and continuum models} Comparison between the simulation results for the discrete model (solid lines) and numerical solutions of the generalised PKS model~\eqref{eq:KellerSegelSystem}-\eqref{eq:exponential_chemofun} (dotted lines). The top panels display the cell density (blue lines) and the bottom panels display the concentration of chemoattractant (red lines) at the end of numerical simulations for progressively smaller cell numbers and critical cell densities, that is, progressively lower values of the parameter $A^0$ in~\eqref{eq:ICBsmall} -- \emph{i.e.} {\bf a.} $A^0 = 10^{5}$, {\bf b.} $A^0 = 10^{4}$, {\bf c.} $A^0 = 10^{3}$, {\bf d.} $A^0 = 10^{2}$. The results from the discrete model correspond to the average over five realisations of the underlying biased random walk. The cell density and the concentration of chemoattractant resulting from each realisation are plotted in pale blue and magenta, respectively. A complete description of the set-up of numerical simulations and the numerical methods employed is given in Appendix~B.}
	\label{fig:lowercell}
\end{figure}

\paragraph{Emergence of differences between spatial patterns produced by the discrete and continuum models}
In all cases discussed so far we have observed excellent agreement between the results of numerical simulations of the discrete model and numerical solutions of the corresponding continuum model given by~\eqref{eq:KellerSegelSystem}-\eqref{eq:exponential_chemofun}. However, we expect possible differences between the two models to emerge in the presence of low cell numbers, which may cause a reduction in the quality of the approximations employed in the formal derivation of the continuum model from the discrete model. To investigate this further, we carry out numerical simulations of the discrete and continuum models considering progressively smaller cell numbers and critical cell densities, \emph{i.e.} defining  
\begin{equation}
\label{eq:ICBsmall}
n_i^0 \equiv A^0, \quad u^0(x) \equiv A^0 \quad  \text{and} \quad u_{\rm max} := 2 \ A^0
\end{equation}
and considering progressively lower values of $A^0$. As expected, the numerical results presented in Figure~\ref{fig:lowercell} show that differences between the patterns produced by the discrete model and those produced by the generalised PKS model~\eqref{eq:KellerSegelSystem}-\eqref{eq:exponential_chemofun} emerge for relatively small cell numbers and critical cell densities, \emph{i.e.} when sufficiently small values of $A^0$ in~\eqref{eq:ICBsmall} are considered. 

\section{Conclusions}
\label{sec:fin}
In this paper, we presented a discrete model of chemotaxis with volume-filling effects. We formally showed that a general form of the celebrated Patlak-Keller-Segel model of chemotaxis can be formally derived as the appropriate continuum limit of this discrete model. Additionally, we characterised the family of steady-state solutions of such a generalised PKS model and we studied the conditions for the emergence of spatial patterns via linear stability analysis. Moreover, we carried out numerical simulations of the discrete and continuum models. 

We showed that there is excellent agreement between the simulation results for the discrete model, the numerical solutions of the corresponding PKS model and the linear stability analysis. Furthermore, we provided numerical evidence for the fact that the dynamics of the cell density and the concentration of chemoattractant exhibited by the two models faithfully replicate those of the classical PKS model in a suitable asymptotic regime. Finally, we showed that possible differences between spatial patterns produced by the two models can emerge in the presence of relatively small cell numbers, which reduce the quality of the approximation of the discrete model given by the continuum model. This demonstrates the importance of integrating discrete and continuum approaches when considering cell-chemotaxis models for spatial pattern formation.

Our discrete modelling framework for chemotactic movement, along with the related formal method to derive corresponding continuum models, can be easily extended to incorporate the effects of additional biological phenomena, such as quorum sensing, haptotaxis and mechanically regulated or nutrient-limited growth of the cell population. An additional development of our study would be to compare the results presented here with those obtained from equivalent models defined on irregular lattices, as well as to investigate how our modelling approach could be related to off-lattice discrete models of cell movement~\cite{drasdo2005coarse}.

\paragraph{Fundin}The authors gratefully acknowledge support of the project PICS-CNRS no. 07688. FB acknowledges funding from the European Research Council (ERC, grant agreement No. 740623) and the Universit\'e Franco-Italienne.

\paragraph{Acknowledgments}The authors are grateful to the two anonymous Reviewers for their useful and insightful comments.

\bibliographystyle{plain} 

\bibliography{references}

\appendix
\section{{Formal derivation of the generalised PKS model~\eqref{eq:KellerSegelSystem}-\eqref{eq:exponential_chemofun}}}
\label{app:derive}

\paragraph{Formal derivation of the conservation equation~\eqref{eq:KellerSegelSystem} for the cell density $u$}
In the case where cell movement is governed by the rules described in Section~\ref{sec:hybrid}, the principle of mass conservation gives 
\begin{equation}
\begin{split}
u^{k+1}_{i} = & \, u^{k}_{i} +  \frac{\theta}{2} \ \psi(u^{k}_{i}) \ \left(u^{k}_{i-1} + u^{k}_{i+1}\right) -\frac{\theta}{2}\left(\psi(u^{k}_{i-1})+\psi(u^{k}_{i+1})\right) u^{k}_{i} \\
& + \frac{\eta}{2\overline{c}} \, \psi(u^{k}_{i})  \left( (c^{k}_{i}-c^{k}_{i-1})_{+} \ u^{k}_{i-1} + (c^{k}_{i}-c^{k}_{i+1})_{+} \ u^{k}_{i+1} \right) \\
& - \frac{\eta}{2\overline{c}} \left(\psi(u^{k}_{i-1}) (c^{k}_{i-1}-c^{k}_{i})_{+}+\psi(u^{k}_{i+1}) (c^{k}_{i+1}-c^{k}_{i})_{+}\right) u^{k}_{i} \ .
\end{split}
\label{eq:intera1}
\end{equation}
Using the fact that the following relations hold for $\tau$ and $h$ sufficiently small
\beq
\label{e:discindcon}
t _k \approx t, \quad t _{k+1} \approx t + \tau, \quad x_{i} \approx x, \quad x_{i \pm 1} \approx x \pm h,
\eeq
\beq
\label{e:discucon}
u^k_{i} \approx u(t,x), \quad u^{k+1}_{i} \approx u(t+\tau,x), \quad u^k_{i \pm 1} \approx u(t,x \pm h),
\eeq
\beq
\label{e:discccon}
c^k_{i} \approx c(t,x), \quad c^{k+1}_{i} \approx c(t+\tau,x), \quad c^k_{i \pm 1} \approx c(t,x \pm h),
\eeq
equation~ \eqref{eq:intera1} can be formally rewritten in the approximate form
\begin{equation*}
\begin{split}
u(t+\tau,x) - u(t,x) = & \, \frac{\theta}{2} \ \psi(u(t,x)) \ \Big(u(t,x-h) + u(t,x+h)\Big) \\ 
& -\frac{\theta}{2}\Big(\psi(u(t,x-h))+\psi(u(t,x+h))\Big) u(t,x) \\
& + \frac{\eta}{2\overline{c}} \, u(t,x-h) \, \psi(u(t,x)) (c(t,x) - c(t,x-h))_{+}\\
& + \frac{\eta}{2\overline{c}} \,u(t,x+h) \, \psi(u(t,x)) (c(t,x) - c(t,x+h))_{+} \\
& - \frac{\eta}{2\overline{c}} \, u(t,x) \, \psi(u(t,x-h)) (c(t,x-h) - c(t,x))_{+}\\
& - \frac{\eta}{2\overline{c}} \, u(t,x)\, \psi(u(t,x+h)) (c(t,x+h) - c(t,x))_{+}.
\end{split}
\end{equation*}
Dividing both sides of the above equation by $\tau$ gives
\begin{equation}
\label{e:deriv1}
\begin{split}
\frac{u(t+\tau,x) - u(t,x)}{\tau} = & \, \frac{\theta}{2 \tau} \ \psi(u(t,x)) \ \Big(u(t,x-h) + u(t,x+h)\Big) \\ 
& -\frac{\theta}{2 \tau}\Big(\psi(u(t,x-h))+\psi(u(t,x+h))\Big) u(t,x) \\
& + \frac{\eta}{2 \tau \overline{c}} \, u(t,x-h) \, \psi(u(t,x)) (c(t,x) - c(t,x-h))_{+}\\
& + \frac{\eta}{2 \tau \overline{c}} \,u(t,x+h) \, \psi(u(t,x)) (c(t,x) - c(t,x+h))_{+} \\
& - \frac{\eta}{2 \tau \overline{c}} \, u(t,x) \, \psi(u(t,x-h)) (c(t,x-h) - c(t,x))_{+}\\
& - \frac{\eta}{2 \tau \overline{c}} \, u(t,x)\, \psi(u(t,x+h)) (c(t,x+h) - c(t,x))_{+}.
\end{split}
\end{equation}
If the function $\psi(u)$ is twice continuously differentiable and the function $u(t,x)$ is twice continuously differentiable with respect to the variable $x$, for $h$ sufficiently small we can use the Taylor expansions
$$
u(t,x \pm h) = u \pm h \frac{\partial u}{\partial x} + \frac{h^2}{2} \frac{\partial^2 u}{\partial x^2} + \mathcal{O}(h^3), \quad \psi(u(t,x \pm h)) = \psi \, \pm \ h \frac{\partial \psi}{\partial x} + \frac{h^2}{2} \frac{\partial^2 \psi}{\partial x^2} +\mathcal{O}(h^3),  
$$
where 
$$
\psi \equiv \psi(u), \quad \frac{\partial \psi}{\partial x} = \psi'(u) \dfrac{\partial u}{\partial x}, \quad \frac{\partial^2 \psi}{\partial x^2} = \psi''(u) \ \left(\frac{\partial u}{\partial x} \right)^2 + \psi'(u) \ \frac{\partial^2 u}{\partial x^2}, \quad u \equiv u(t,x).
$$
Substituting into~\eqref{e:deriv1} and using the elementary property $(a)_+ - (-a)_+ = a$ for $a \in \mathbb{R}$, after a little algebra we find
\begin{equation*}
\begin{split}
\frac{u(t+\tau,x) - u(t,x)}{\tau}  = & \,  \frac{\theta h^{2}}{2 \tau}\psi\ \frac{\partial^2 u}{\partial x^2}-\frac{\theta h^{2}}{2 \tau} \frac{\partial^2 \psi}{\partial x^2} \ u \\
& +\frac{\eta}{2 \tau \overline{c}} \ \psi \ \Big(2c(t,x) - c(t,x-h) - c(t,x+h) \Big) \ u\\
& -\frac{\eta h}{2 \tau \overline{c}} \ \psi \ \Big(\left(c(t,x)-c(t,x-h)\right)_{+} -\left(c(t,x)-c(t,x+h)\right)_{+}\Big) \ \frac{\partial u}{\partial x} \\
& +\frac{\eta h}{2 \tau \overline{c}} \ \frac{\partial \psi}{\partial x} \ \Big(\left(c(t,x-h)-c(t,x)\right)_{+} - \left(c(t,x+h)-c(t,x)\right)_{+} \Big) \ u\\
& +\frac{\eta h^2}{4 \tau \overline{c}} \ \psi \ \Big(\left(c(t,x)-c(t,x-h)\right)_{+} + \left(c(t,x)-c(t,x+h)\right)_{+} \Big) \ \frac{\partial^2 u}{\partial x^2} \\
& -\frac{\eta h^2}{4 \tau \overline{c}} \ \frac{\partial^2 \psi}{\partial x^2} \ \Big(\left(c(t,x-h)-c(t,x)\right)_{+} + \left(c(t,x+h)-c(t,x)\right)_{+}\Big) \ u \ + h.o.t. \, ,
\end{split}
\end{equation*}
which can be rewritten as
\begin{equation*}
\begin{split}
\frac{u(t+\tau,x) - u(t,x)}{\tau} = & \,  \frac{\theta h^{2}}{2 \tau}\ \left(\psi \frac{\partial^2 u}{\partial x^2} - \frac{\partial^2 \psi}{\partial x^2} \ u\right) \\
& -\frac{\eta h^2}{2 \tau \overline{c}} \ \psi \ \left(\frac{c(t,x-h) + c(t,x+h) - 2 c(t,x)}{h^2} \right) \ u\\
& -\frac{\eta h^2}{2 \tau \overline{c}} \ \psi \ \left(\left(\frac{c(t,x)-c(t,x-h)}{h}\right)_{+} -\left(\frac{c(t,x)-c(t,x+h)}{h}\right)_{+}\right) \ \frac{\partial u}{\partial x} \\
& +\frac{\eta h^2}{2 \tau \overline{c}} \ \frac{\partial \psi}{\partial x} \ \left(\left(\frac{c(t,x-h)-c(t,x)}{h}\right)_{+} - \left(\frac{c(t,x+h)-c(t,x)}{h}\right)_{+} \right) \ u\\
& + \frac{\eta h^3}{4 \tau \overline{c}} \ \psi \ \left(\left(\frac{c(t,x)-c(t,x-h)}{h}\right)_{+} + \left(\frac{c(t,x)-c(t,x+h)}{h}\right)_{+} \right) \ \frac{\partial^2 u}{\partial x^2} \\
& - \frac{\eta h^3}{4 \tau \overline{c}} \ \frac{\partial^2 \psi}{\partial x^2} \ \left(\left(\frac{c(t,x-h)-c(t,x)}{h}\right)_{+} + \left(\frac{c(t,x+h)-c(t,x)}{h}\right)_{+}\right) \ u \\
& + h.o.t. \, .
\end{split}
\end{equation*}
If, in addition, the function $u(t,x)$ is continuously differentiable with respect to the variable $t$ and the function $c(t,x)$ is twice continuously differentiable with respect to the variable $x$, letting $\tau \to 0$ and $h \to 0$ in such a way that
\begin{equation}
\label{ass:quottozero1dapp}
\frac{\eta h^2}{2\tau \overline{c}}\rightarrow \chi \in \mathbb{R}^+_* \quad \text{and} \quad \frac{\theta h^2}{2\tau}\rightarrow \beta_{u} \in \mathbb{R}^+_* \quad \text{as} \quad \tau, h \to 0,
\end{equation}
we have 
\begin{equation*}
\frac{\eta h^3}{4 \tau \overline{c}} \psi \left(\left(\frac{c(t,x)-c(t,x-h)}{h}\right)_{+} + \left(\frac{c(t,x)-c(t,x+h)}{h}\right)_{+} \right) \frac{\partial^2 u}{\partial x^2} = \mathcal{O}(h), \quad \text{as } \tau, h \to 0,
\end{equation*}
\begin{equation*}
\frac{\eta h^3}{4 \tau \overline{c}} \frac{\partial^2 \psi}{\partial x^2} \left(\left(\frac{c(t,x-h)-c(t,x)}{h}\right)_{+} + \left(\frac{c(t,x+h)-c(t,x)}{h}\right)_{+}\right) u = \mathcal{O}(h), \quad \text{as } \tau, h \to 0,
\end{equation*}
and from the latter equation we formally obtain
\begin{equation*}
\begin{split}
\frac{\partial u}{\partial t} = & \,  \beta_{u}\ \Big(\psi \frac{\partial^2 u}{\partial x^2}  -u \ \frac{\partial^2 \psi}{\partial x^2} \Big) \\
& - \chi \ \psi \ u \ \frac{\partial^2 c}{\partial x^2} - \chi \ \psi \ \frac{\partial u}{\partial x} \ \left(\left(\frac{\partial c}{\partial x}\right)_{+} -\left(-\frac{\partial c}{\partial x}\right)_{+}\right) + \chi \ \frac{\partial \psi}{\partial x} \ u \ \left(\left(-\frac{\partial c}{\partial x}\right)_{+} - \left(\frac{\partial c}{\partial x}\right)_{+} \right).
\end{split}
\end{equation*}
Using again the elementary property $(a)_+ - (-a)_+ = a$ for $a \in \mathbb{R}$, we find
\begin{equation}
\label{eq:finalu}
\frac{\partial u}{\partial t} =  \beta_{u}\ \Big(\psi \frac{\partial^2 u}{\partial x^2}  -u \ \frac{\partial^2 \psi}{\partial x^2} \Big) -\chi \frac{\partial }{\partial x}\Big( \psi \ u\ \frac{\partial c}{\partial x}\Big),
\end{equation}
where $\psi \equiv \psi(u)$, $u \equiv u(t,x)$ and $c \equiv c(t,x)$. Since
\begin{equation*}
\psi(u) \frac{\partial^2 u}{\partial x^2}  - u \frac{\partial^2 \psi(u)}{\partial x^2} = \frac{\partial }{\partial x} \left[ \left( \psi(u)-u\ \psi^{\prime}(u)\right) \frac{\partial u}{\partial x}\right],
\end{equation*}
under definition~\eqref{ass:volfil} of the function $D(u)$ the differential equation~\eqref{eq:finalu} can be rewritten as
\begin{equation*}
\frac{\partial u}{\partial t} - \frac{\partial }{\partial x} \left(\beta_u \, D(u) \, \frac{\partial u}{\partial x} -  \chi \, \psi(u) \, u \, \frac{\partial c}{\partial x}   \right) = 0,
\end{equation*}
which is the conservation equation~\eqref{eq:KellerSegelSystem} for the cell density $u$ complemented with~\eqref{ass:volfil} and~\eqref{eq:exponential_chemofun}, and posed on $\mathbb{R}^+_* \times \mathbb{R}$.

\paragraph{Formal derivation of the balance equation~\eqref{eq:KellerSegelSystem} for the chemoattractant concentration $c$} For the one-dimensional case considered in Section~\ref{sec:hybrid} we have
$$
(\mathcal{L} \ c^k)_i = \frac{c^{k}_{i+1} - 2 c^{k}_{i} + c^{k}_{i-1}}{h^2}.
$$
Hence, if $\tau$ and $h$ are sufficiently small so that relations~\eqref{e:discindcon}-\eqref{e:discccon} hold, the difference equation~\eqref{e:cdiscrete} can be formally written in the approximate form
$$
\frac{c(t+\tau,x) - c(t,x)}{\tau} = \left(\beta_c \frac{c(t,x-h) + c(t,x+h) - 2 c(t,x)}{h^2} + \alpha \ u(t,x) - \kappa \ c(t,x) \right).
$$
If the function $c(t,x)$ is continuously differentiable with respect to the variable $t$ and twice continuously differentiable with respect to the variable $x$, letting $\tau \to 0$ and $h \to 0$ in the above equation formally gives
$$
\frac{\partial c}{\partial t} - \beta_c \, \frac{\partial^2 c}{\partial x^2} = \alpha \, u - \kappa \, c,
$$
which is the balance equation~\eqref{eq:KellerSegelSystem} for the chemoattractant concentration $c$ posed on $\mathbb{R}^+_* \times \mathbb{R}$.

\section{Details of numerical simulations in 1D}
\label{numerics1D}

\subsection{Details of numerical simulations of the discrete model in 1D}
\label{numerics1DIB}
We use a uniform discretisation of the interval $\Omega := (0,1)$ that consists of $N=100$ points as the spatial domain (\emph{i.e.} the grid-step is $h\approx1 \times 10^{-2}$) and we choose the time-step $\tau=1\times10^{-2}$. Numerical simulations are performed for $5 \times 10^5$ time-steps (\emph{i.e.} the final time of simulations is $t = 500$). 

\begin{figure}[!h]
	\centering
	\includegraphics[width=\textwidth]{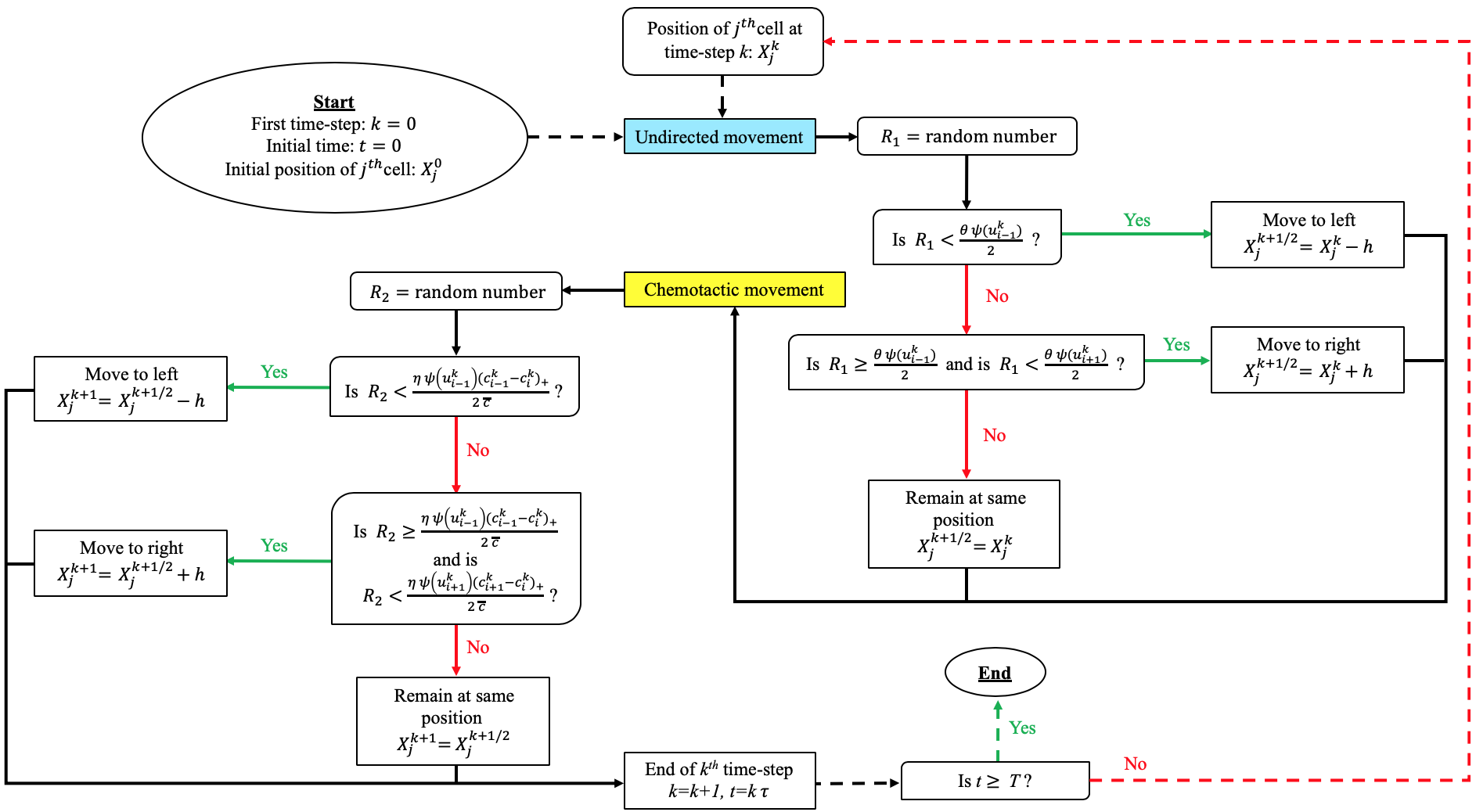}
	\caption{Flowchart illustrating the computational procedure followed to update the positions of every cell in 1D.}
	\label{fig:flowchart}
\end{figure}

\paragraph{Computational procedure} At each time-step, we follow the computational procedure illustrated by the flowchart in Figure~\ref{fig:flowchart} to update the positions of the single cells. Zero-flux boundary conditions are implemented by letting the attempted move of a cell be aborted if it requires moving out of the spatial domain. Numerical simulations are performed in {\sc Matlab} and the random numbers mentioned in Figure~\ref{fig:flowchart} are all real numbers drawn from the standard uniform distribution on the interval $(0,1)$ using the built-in function {\sc rand}. At each time-step, the positions of all the cells are updated first and then the cell density at every lattice site is computed via~\eqref{e:n} and inserted into~\eqref{e:cdiscrete} in order to update the concentration of the chemoattractant.

\paragraph{Initial conditions} 
The numerical results of Figures~\ref{fig:basecase},~\ref{fig:chi} and~\ref{fig:infinity} refer to the case where the cells are initially uniformly distributed and the corresponding cell density is
\begin{equation*} 
u^0_{i} = \frac{A}{2} \; \text{ for } \;  i=1,\ldots,N, \quad A=2\times10^{6}.
\end{equation*}
Moreover, the numerical results of Figure~\ref{fig:initial} refer to the case where
\begin{equation*}
u^0_{i} = \frac{A}{2} B \; \text{ for } \;  i=1,\ldots,N, \quad A=4\times10^{5},
\end{equation*}
and the value of $B$ is varied as described in the caption of the figure. {Finally,  the numerical results of Figure~\ref{fig:lowercell} refer to the case where the initial cell density is defined according to~\eqref{eq:ICBsmall} and the value of $A^0$ is varied as described in the caption of the figure.}

{For all cases, the initial concentration of the chemoattractant is defined as an oscillating perturbation around the initial cell density, which is uniform, in order to drive pattern formation, \emph{i.e.} 
\begin{equation*}
c^0_{i}= u^0_{i} \left(1 + 0.1 \ \cos(10\ i \ h)\ \sin(10\ i \ h)\right), \quad  i=1,\ldots,N.
\end{equation*}
An alternative initial condition for the concentration of chemoattractant could be a small random perturbation around $u^0_i$. Preliminary numerical simulations showed that the same spatial patterns are formed in both cases.}

\paragraph{Parameter values} 
For all cases, the following parameter values are used
\begin{equation*}
\beta_{c} = 2.5 \times 10^{-3}, \quad \alpha=1, \quad \kappa=1, \quad \theta=1.225 \times 10^{-1},
\end{equation*}
The numerical results of Figures~\ref{fig:basecase}~and~\ref{fig:chi} refer to the case where $u_{\text{max}}=2\times 10^{6}$, while the numerical results of Figure~\ref{fig:initial} refer to the case where $u_{\text{max}}=4\times 10^{5}$. For the numerical results of Figure~\ref{fig:infinity} the value of $u_{\text{max}}$ is varied as described in the caption of the figure and {for Figure~\ref{fig:lowercell} the value of $u_{\text{max}}$ is defined according to~\eqref{eq:ICBsmall}, with the value of $A^0$ being varied as described in the caption of the figure.} In all cases, the value of $\bar{c}$ is defined via~\eqref{e:defbarc} with $\zeta=1$. The numerical results presented in Figures~\ref{fig:basecase},~\ref{fig:initial}~and~~\ref{fig:lowercell}, refer to the case where $\eta = 2.4502$, while for the numerical results of Figure~\ref{fig:chi} the value of $\eta$ is varied as described in the caption of the figure. Furthermore, for the results presented in Figure~\ref{fig:infinity}, the value of $\eta$ is varied according to the value of $u_{\text{max}}$ so that their quotient remains constant and equal to $1.225\times10^{-6}$.

\subsection{Details of numerical simulations of the generalised PKS model~\eqref{eq:KellerSegelSystem}-\eqref{eq:exponential_chemofun} in 1D}
We select a uniform discretisation consisting of $\mathcal{N}=100$ points of the interval $\Omega := (0,1)$ as the computational domain of the independent variable $x$ (\emph{i.e.} $x_j = j \, \Delta x$ with $\Delta x \approx1 \times 10^{-2}$ and $j=1, \ldots, \mathcal{N}$). Moreover, we assume $t \in (0,500]$ and we discretise the interval $(0,500]$ with the uniform step $\Delta t = 1 \times 10^{-2}$.

\paragraph{Numerical methods} The method for constructing numerical solutions of the generalised PKS model~\eqref{eq:KellerSegelSystem}-\eqref{eq:exponential_chemofun} (\emph{i.e.} the system of coupled parabolic equations~\eqref{modexpli}) is based on a finite difference scheme whereby the discretised dependent variables are
\begin{equation*}
	u^n_j \approx u(t_n, x_j) \quad \text{and} \quad c^n_j \approx c(t_n, x_j).
\end{equation*}
We solve numerically equation~\eqref{eq:KellerSegelSystem} for $c$ using an implicit Euler method, that is, 
\begin{equation*}
	\frac{c_j^{n+1}-c_j^{n}}{\Delta t} = \beta_c \frac{c^{n+1}_{j+1}-2c_j^{n+1} + c_{j-1}^{n+1}}{\left(\Delta x \right)^2} + \kappa c_j^{n+1} - \alpha u_j^{n}, \quad j = 1, \dots, \mathcal{N},
\end{equation*}
and impose zero-flux boundary conditions by letting
\begin{equation*}
c^{n+1}_{0} = c^{n+1}_{1} \quad \text{and} \quad  c^{n+1}_{\mathcal{N}+1} = c^{n+1}_{\mathcal{N}}.
\end{equation*}
Moreover, we solve numerically equation~\eqref{eq:KellerSegelSystem} for $u$ using the following implicit scheme
\begin{equation*}
\frac{u^{n+1}_j - u^n_j}{\Delta t} = \frac{F_{j+\frac{1}{2}}^{n+1} - F_{j-\frac{1}{2}}^{n+1}}{\Delta x}, \quad j = 1, \dots, \mathcal{N},
\end{equation*}
where 
\begin{equation*}
\begin{split}
F_{j + \frac{1}{2}}^{n+1} &:= \beta_u D\left(u^n_{j+\frac{1}{2}}\right)\frac{u^{n+1}_{j+1}-u^{n+1}_{j}}{\Delta x}  - b^{n,+}_{j+\frac{1}{2}} u^{n+1}_{j} \psi\left(u^{n+1}_{j+1}\right)\\
						& \quad  + b^{n,-}_{j+\frac{1}{2}} u^{n+1}_{j+1} \psi\left(u^{n+1}_{j}\right), \quad j = 1, \dots, \mathcal{N}-1,
\end{split}
\end{equation*}
with
\begin{equation*}
u^n_{j+\frac{1}{2}} := \frac{u^n_{j+1}+u^n_{j}}{2},
\end{equation*}
and
\begin{equation*}
b^{n}_{j+\frac{1}{2}} := \chi \frac{c^n_{j+1}-c^n_{j}}{\Delta x}, \quad b^{n,+}_{j+\frac{1}{2}} = \max\left(0,b^{n}_{j+\frac{1}{2}} \right), \quad b^{n,-}_{j+\frac{1}{2}} = \max\left(0,-b^{n}_{j+\frac{1}{2}} \right).
\end{equation*}
The discrete flux $F_{j-\frac{1}{2}}^{n+1}$ for $j = 2, \dots, \mathcal{N}$ is defined in an analogous way, and we impose zero-flux boundary conditions by using the definitions
\begin{equation*}
F_{1 - \frac{1}{2}}^{n+1} := 0 \quad \text{and} \quad F_{\mathcal{N} + \frac{1}{2}}^{n+1} := 0.
\end{equation*}
Analogous schemes are used to solve numerically the classical PKS model~\eqref{eq:KellerSegelSystemOri}. All numerical computations are performed in {\sc Matlab}. 

\paragraph{Initial conditions and parameter values} 
In agreement with the set-up of numerical simulations of the discrete model, the numerical results of Figures~\ref{fig:basecase},~\ref{fig:chi} and~\ref{fig:infinity} refer to the case where
\begin{equation*}
u(0,x) \equiv \frac{A}{2}, \quad A=2\times10^{6}.
\end{equation*}
Moreover, the numerical results of Figure~\ref{fig:initial} refer to the case where
\begin{equation*}
u(0,x) \equiv \frac{A}{2} B, \quad A=4\times10^{5},
\end{equation*}
and the value of $B$ is varied as described in the caption of the figure. {Finally,  the numerical results of Figure~\ref{fig:lowercell} refer to the case where the initial cell density is defined according to~\eqref{eq:ICBsmall} and the value of $A^0$ is varied as described in the caption of the figure.} In all cases, 
\begin{equation*}
c(0,x)= u(0,x)+0.1\left(u(0,x) \ \cos(10\ x)\ \sin(10\ x)\right).
\end{equation*}
\paragraph{Parameter values} 
In agreement with the set-up of numerical simulations of the discrete model, for all cases the following parameter values are used
\begin{equation*}
\beta_{c} = 2.5 \times 10^{-3}, \quad \alpha=1, \quad \kappa=1.
\end{equation*}
The numerical results of Figures~\ref{fig:basecase}~and~\ref{fig:chi} refer to the case where $u_{\text{max}}=2\times 10^{6}$, while the numerical results of Figure~\ref{fig:initial} refer to the case where $u_{\text{max}}=4\times 10^{5}$. For the numerical results of Figure~\ref{fig:infinity} the value of $u_{\text{max}}$ is varied as described in the caption of the figure and {for Figure~\ref{fig:lowercell} the value of $u_{\text{max}}$ is defined according to~\eqref{eq:ICBsmall}, with the value of $A^0$ being varied as described in the caption of the figure.} In all cases, given the values of the parameters chosen to carry out numerical simulations of the discrete model (see Appendix~\ref{numerics1D}\ref{numerics1DIB}), the following definitions are used
\begin{equation*}
\chi := \frac{\eta h^2}{2\tau \overline{c}} \quad \text{and} \quad \beta_u:= \frac{\theta h^2}{2\tau},
\end{equation*}
so that conditions~\eqref{ass:quottozero1d} are met.

\section{Details of numerical simulations in 2D}
\label{numerics2D}
\subsection{Details of numerical simulations of the discrete model in 2D}
\label{numerics2DIB}
We use a uniform discretisation of the square $\Omega := (0,1) \times (0,1)$ that consists of $N^2=2601$ points as the spatial domain (\emph{i.e.} the grid-step is $h\approx1.9 \times 10^{-2}$) and we choose the time-step $\tau=1\times10^{-4}$. Numerical simulations are performed for $5 \times 10^4$ time-steps (\emph{i.e.} the final time of simulations is $t = 5$) for the numerical results of Figure~\ref{fig:infinity2D}, whereas in all the other cases simulations are performed for $15 \times 10^4$ time-steps (\emph{i.e.} the final time of simulations is $t = 15$).   

\paragraph{Computational procedure} At each time-step, the positions of the single cells are updated following a procedure analogous to that employed in 1D (\emph{cf.} Figure~\ref{fig:flowchart}), with the only difference being that the cells are allowed to move up and down as well. Moreover, the concentration of the chemoattractant is updated through the two-dimensional analogue of~\eqref{e:cdiscrete}, where the operator $\mathcal{L}$ is defined as the finite-difference Laplacian on a two-dimensional regular grid of step $h$ and the cell density is computed via the two-dimensional analogue of~\eqref{e:n}.  

\paragraph{Initial conditions} 
The numerical results of Figures~\ref{fig:basecase2D} and~\ref{fig:infinity2D} refer to the case where the cells are initially uniformly distributed and the corresponding cell density is
\begin{equation*}
u^0_{ij} = \frac{A}{2} \; \text{ for } \; i,j=1,\ldots,N, \quad A=1\times10^{7},
\end{equation*}
while the numerical results of Figure~\ref{fig:initial2D} refer to the case where
\begin{equation*}
u^0_{ij} = \frac{A}{2} B \; \text{ for } \; i,j=1,\ldots,N, \quad A=1\times10^{7},
\end{equation*}
and the value of $B$ is varied as described in the caption of the figure. For the numerical results of Figures~\ref{fig:basecase2D},~\ref{fig:initial2D}~and~\ref{fig:infinity2D}, the initial concentration of chemoattractant is defined as
\begin{equation*}
c^0_{ij}= 200 \sum_{p=1}^4 \exp\left[-200 (i \, h-x^{*}_{1p})^{2} -200 (j \, h-x^{*}_{2p})^{2} \right], \quad i,j=1,\ldots,N.
\end{equation*}
Furthermore, in all cases, $(x^{*}_{11}, x^{*}_{21})=(0.26, 0.74)$, $(x^{*}_{12}, x^{*}_{22})=(0.26,0.26)$, $(x^{*}_{13}, x^{*}_{23})=(0.74,0.74)$ and $(x^{*}_{14}, x^{*}_{24})=(0.74,0.26)$.

\paragraph{Parameter values} 
For all cases, the following parameter values are used
\begin{equation*}
\beta_{c} = 2.5 \times 10^{-3}, \quad \alpha=1, \quad \kappa=1.
\end{equation*}
Moreover, $\theta=2.5 \times 10^{-2}$ for the results presented in Figures~\ref{fig:basecase2D}~and~\ref{fig:initial2D}, while $\theta=0.125$ for the results presented in Figure~\ref{fig:infinity2D}. The numerical results of Figures~\ref{fig:basecase2D}~and~\ref{fig:initial2D} we choose $u_{\text{max}}=10^{7}$, and for the results of Figure~\ref{fig:infinity2D} the value of $u_{\text{max}}$ is varied as described in the caption of the figure. In all cases, the value of $\bar{c}$ is defined via~\eqref{e:defbarc} with $\zeta=1$. The numerical results presented in Figures~\ref{fig:basecase2D} and~\ref{fig:initial2D} refer to the case where $\eta = 2.4502$, while for the numerical results presented in Figure~\ref{fig:infinity}, the value of $\eta$ is varied according to the value of $u_{\text{max}}$ so that their quotient remains constant and equal to $2.4502\times10^{-7}$.

\subsection{Details of numerical simulations of the generalised PKS model~\eqref{eq:KellerSegelSystem}-\eqref{eq:exponential_chemofun} in 2D}
We select a uniform discretisation consisting of $\mathcal{N}^2 = 2601$ points of the square $\Omega := (0,1) \times (0,1)$ as the computational domain of the independent variable $x \equiv (x_1, x_2)$ (\emph{i.e.} $(x_{1 i}, x_{2 j}) = \left(i \, \Delta x, j \, \Delta x\right)$ with $\Delta x \approx 1.9 \times 10^{-2}$ and $i,j=1, \ldots, \mathcal{N}$). Moreover, we assume $t \in (0,T]$ with $T=5$ for the numerical solutions of Figure~\ref{fig:infinity2D},  whereas in all the other cases $T=15$. The interval $(0,T]$ is discretised with the uniform step $\Delta t =1\times10^{-4}$.

\paragraph{Numerical methods} The method for constructing numerical solutions of the generalised PKS model~\eqref{eq:KellerSegelSystem}-\eqref{eq:exponential_chemofun} (\emph{i.e.} the system of coupled parabolic equations~\eqref{modexpli}) is based on a finite difference scheme whereby the discretised dependent variables are
\begin{equation*}
	u^n_{i,j} := u(t_n, x_{1 i}, x_{2 j}) \quad \text{and} \quad c^n_{i,j} := c(t_n, x_{1 i}, x_{2 j}).
\end{equation*}
We solve numerically equation~\eqref{eq:KellerSegelSystem} for $c$ using an implicit Euler method, that is, 
\begin{equation*}
\begin{split}
\frac{c^{n+1}_{i,j}-c^{n}_{i,j}}{\Delta t} &= \beta_c \frac{c^{n}_{i+1,j}-2c^{n}_{i,j} + c^{n}_{i-1,j}}{\left(\Delta x\right)^2} +  \beta_c \frac{c^{n}_{i,j+1}-2c^{n}_{i,j}+c^{n}_{i,j-1}}{\left(\Delta x\right)^2} \\
						&\quad + \alpha u^{n}_{i,j} - \kappa c^n_{i,j}, \quad  i,j = 1, \dots, \mathcal{N},
\end{split}
\end{equation*}
and impose zero-flux boundary conditions by letting
\begin{equation*}
	\begin{split}
	c^{n+1}_{0,j} = c^{n+1}_{1,j}, \quad c^{n+1}_{\mathcal{N}+1,j} = c^{n+1}_{\mathcal{N},j}, \quad j = 1, \dots, \mathcal{N},\\
	c^{n+1}_{i,0} = c^{n+1}_{i,1}, \quad c^{n+1}_{i,\mathcal{N}+1} = c^{n+1}_{i,\mathcal{N}}, \quad i = 1, \dots, \mathcal{N}.
	\end{split}
\end{equation*}
Moreover, we solve numerically equation~\eqref{eq:KellerSegelSystem} for $u$ using the explicit scheme
\begin{equation*}
	\frac{u^{n+1}_{i,j} - u^n_{i,j}}{\Delta t} = \frac{F^{n}_{i+\frac{1}{2},j}- F^{n}_{i-\frac{1}{2},j}}{\Delta x} + \frac{F^{n}_{i,j+\frac{1}{2}}-F^{n}_{i,j-\frac{1}{2}}}{\Delta x}, \quad i, j = 1, \dots, \mathcal{N},
\end{equation*}
where 
\begin{equation*}
\begin{split}
F^{n}_{i+\frac{1}{2},j} &: = \beta_u D\left(u^n_{i+\frac{1}{2},j}\right) \frac{u^n_{i+1,j}-u^n_{i,j}}{\Delta x} - b^{n,+}_{i+\frac{1}{2},j} u^n_{i,j} \psi(u^n_{i+1,j})\\
						& \quad  +b^{n,-}_{i+\frac{1}{2},j} u^n_{i+1,j} \psi(u^n_{i,j}), \quad i = 1,\dots,\mathcal{N}-1, \,\, j = 1, \dots, \mathcal{N},\\
F^{n}_{i,j+\frac{1}{2}} &: = \beta_u D\left(u^n_{i,j+\frac{1}{2}}\right) \frac{u^n_{i,j+1}-u^n_{i,j}}{\Delta x} - b^{n,+}_{i,j+\frac{1}{2}} u^n_{i,j} \psi(u^n_{i,j+1})\\
& \quad  +b^{n,-}_{i,j+\frac{1}{2}} u^n_{i,j+1} \psi(u^n_{i,j}), \quad i = 1,\dots,\mathcal{N}, \,\, j = 1, \dots, \mathcal{N}-1,
\end{split}
\end{equation*}
with
\begin{equation*}
u^n_{i+\frac{1}{2},j} := \frac{u^n_{i+1,j}+u^n_{i,j}}{2}, \quad u^n_{i,j+\frac{1}{2}} := \frac{u^n_{i,j+1}+u^n_{i,j}}{2},
\end{equation*}
\begin{equation*}
b^{n}_{i+\frac{1}{2},j} : = \chi \frac{c^n_{i+1,j}-c^n_{i,j}}{\Delta x}, \quad 	b^{n,+}_{i+\frac{1}{2},j} = \max \left(0, b^{n}_{i+\frac{1}{2},j} \right), \quad b^{n,-}_{i+\frac{1}{2},j} = \max \left(0, -b^{n}_{i+\frac{1}{2},j} \right),
\end{equation*}
and
\begin{equation*}
b^{n}_{i,j+\frac{1}{2}} : = \chi \frac{c^n_{i,j+1}-c^n_{i,j}}{\Delta x}, \quad b^{n,+}_{i,j+\frac{1}{2}} = \max \left(0, b^{n}_{i,j+\frac{1}{2}} \right), \quad b^{n,-}_{i,j+\frac{1}{2}} = \max \left(0, -b^{n}_{i,j+\frac{1}{2}} \right).
\end{equation*}
The discrete fluxes $F^{n}_{i-\frac{1}{2},j}$ for $i = 2,\dots,\mathcal{N}$, $j = 1, \dots, \mathcal{N}$ and $F^{n}_{i,j-\frac{1}{2}}$ for $i = 1,\dots,\mathcal{N}$, $j = 2, \dots, \mathcal{N}$ are defined in analogous ways, and we impose zero-flux boundary conditions by using the definitions
\begin{equation*}
\begin{split}
F^{n}_{1-\frac{1}{2},j} := 0, \quad F^{n}_{\mathcal{N}+\frac{1}{2},j} := 0, \quad j = 1, \dots, \mathcal{N},\\
F^{n}_{i,1-\frac{1}{2}} := 0,  \quad F^{n}_{i,\mathcal{N}+\frac{1}{2}} := 0, \quad i = 1, \dots, \mathcal{N}.
\end{split}
\end{equation*}
Notice that, in contrast to the one-dimensional case, here we employ a fully explicit scheme to avoid Newton sub-iterations that could be computationally expensive. Analogous schemes are used to solve numerically the classical PKS model~\eqref{eq:KellerSegelSystemOri}. All numerical computations are performed in {\sc Matlab}. 

\paragraph{Initial conditions} 

In agreement with the set-up of numerical simulations of the discrete model, the numerical results of Figures~\ref{fig:basecase2D} and~\ref{fig:infinity2D} refer to the case where
\begin{equation*}
u(0,x_1,x_2) \equiv \frac{A}{2}, \quad A=1\times10^{7},
\end{equation*}
while the numerical results of Figure~\ref{fig:initial2D} refer to the case where
\begin{equation*}
u(0,x_1,x_2) \equiv \frac{A}{2} B, \quad A=1\times10^{7},
\end{equation*}
and the value of $B$ is varied as described in the caption of the figure. For all cases, the initial concentration of chemoattractant is
\begin{equation*}
c(0,x_1,x_2) = 200 \sum_{p=1}^4 \exp\left[-200 (x_{1}-x^{*}_{1p})^{2} -200 (x_{2}-x^{*}_{2p})^{2} \right],
\end{equation*}
with $(x^{*}_{11}, x^{*}_{21})=(0.26, 0.74)$, $(x^{*}_{12}, x^{*}_{22})=(0.26,0.26)$, $(x^{*}_{13}, x^{*}_{23})=(0.74,0.74)$ and $(x^{*}_{14}, x^{*}_{24})=(0.74,0.26)$.

\paragraph{Parameter values} 
For all cases, the following parameter values are used
\begin{equation*}
\beta_{c} = 0.0025, \quad \alpha=1, \quad \kappa=1.
\end{equation*}
The numerical results of Figures~\ref{fig:basecase2D} and~\ref{fig:initial2D} refer to the case where $u_{\text{max}}=1\times 10^{7}$, while for the numerical results of Figure~\ref{fig:infinity2D} the value of $u_{\text{max}}$ is varied as described in the caption of the figure. In all cases, given the values of the parameters chosen to carry out numerical simulations of the discrete model (see Appendix~\ref{numerics2D}\ref{numerics2DIB}), the following definitions are used
\begin{equation*}
\chi := \frac{\eta h^2}{4\tau \overline{c}} \quad \text{and} \quad \beta_u:= \frac{\theta h^2}{4\tau},
\end{equation*}
so that conditions~\eqref{ass:quottozero2d} are met.

\end{document}